\documentclass[11pt]{amsart}

\newcommand{\FF}{\mathbb{F}}

\newcommand{\NN}{\mathbb{N}}

\newcommand{\ZZ}{\mathbb{Z}}

\newcommand{\ccD}{\mathcal{D}}

\newcommand{\cH}{\mathcal{H}}

\newcommand{\cM}{\mathcal{M}}

\newcommand{\cO}{\mathcal{O}}
\newcommand{\cP}{\mathcal{P}}

\newcommand{\cT}{\mathcal{T}}

\newcommand{\cV}{\mathcal{V}}

\newcommand{\cX}{\mathcal{X}}

\newcommand{\fA}{\mathfrak{A}}

\newcommand{\fc}{\mathfrak{c}}

\newcommand{\Car}{\mathfrak{Car}}

\newcommand{\fg}{\mathfrak{g}}
\newcommand{\fh}{\mathfrak{h}}
\newcommand{\fl}{\mathfrak{l}}
\newcommand{\fm}{\mathfrak{m}}

\newcommand{\fn}{\mathfrak{n}}

\newcommand{\fs}{\mathfrak{s}}

\newcommand{\ft}{\mathfrak{t}}
\newcommand{\fT}{\mathfrak{T}}
\newcommand{\Tor}{\mathfrak{Tor}}

\newcommand{\sfd}{{\textsf{d}}}

\newcommand{\dact}{\boldsymbol{.}}
\newcommand{\lra}{\longrightarrow}

\DeclareMathOperator{\ad}{ad}

\DeclareMathOperator{\Aut}{Aut}
\DeclareMathOperator{\Cent}{Cent}
\DeclareMathOperator{\ch}{char}
\DeclareMathOperator{\Div}{div}
\DeclareMathOperator{\End}{End}
\DeclareMathOperator{\GL}{GL}

\DeclareMathOperator{\Hom}{Hom}
\DeclareMathOperator{\im}{im}
\DeclareMathOperator{\id}{id}

\DeclareMathOperator{\Lie}{Lie}

\DeclareMathOperator{\Nor}{Nor}
\DeclareMathOperator{\Rac}{Rac}
\DeclareMathOperator{\Rad}{Rad}
\DeclareMathOperator{\rk}{rk}
\DeclareMathOperator{\res}{res}

\DeclareMathOperator{\SL}{SL}

\DeclareMathOperator{\Sp}{Sp}

\DeclareMathOperator{\Spec}{Spec}

\usepackage{pictex}
\usepackage{amscd}
\usepackage{latexsym}
\usepackage{amssymb}
\usepackage{eucal}
\usepackage{verbatim}

\numberwithin{equation}{section}

\newtheorem{Theorem}{Theorem}[section]
\newtheorem{Lemma}[Theorem]{Lemma}
\newtheorem{Corollary}[Theorem]{Corollary}
\newtheorem{Proposition}[Theorem]{Proposition}

\newtheorem{Thm}{Theorem}[subsection]
\newtheorem{Lem}[Thm]{Lemma}
\newtheorem{Prop}[Thm]{Proposition}
\newtheorem{Cor}[Thm]{Corollary}

\newtheorem*{thm*}{Theorem A}
\newtheorem*{thm**}{Theorem B}

\theoremstyle{remark}
\newtheorem*{Remark}{Remark}
\newtheorem*{Remarks}{Remarks}
\newtheorem*{Definition}{Definition}

\numberwithin{equation}{section}

\begin{document}

\title[Weyl Groups]{Weyl groups for non-classical Restricted Lie algebras\\ and the Chevalley Restriction Theorem}

\author[J.-M. Bois, R. Farnsteiner \lowercase{and} B. Shu]{Jean-Marie Bois, Rolf Farnsteiner \lowercase{and} Bin Shu}

\address[Jean-Marie Bois]{Mathematisches Seminar, Christian-Albrechts-Universit\"at zu Kiel, Ludewig-Meyn-Str.~4, 24098 Kiel, Germany}
\email{bois@math.uni-kiel.de}

\address[Rolf Farnsteiner]{Mathematisches Seminar, Christian-Albrechts-Universit\"at zu Kiel, Ludewig-Meyn-Str.~4, 24098 Kiel, Germany}
\email{rolf@math.uni-kiel.de}
\thanks{The first and second named authors were supported by the D.F.G. priority program SPP1388 `Darstellungstheorie'. The third named author was partially supported by NSF of China (Grant No:10871067), PCSIRT. He thanks the DAAD, as well as the Department of Mathematics and the SFB at the University of Bielefeld for the hospitality and support during his visit in 2007--2008.}
\address[Bin Shu]{Department of Mathematics, East China Normal University, Shanghai 200062,  China}
\email{bshu@math.ecnu.edu.cn}

\date{\today}

\makeatletter
\makeatother


\begin{abstract}
Let $(\fg,[p])$ be a finite-dimensional restricted Lie algebra, defined over an algebraically closed field $k$ of characteristic $p>0$. The scheme of tori of maximal dimension of $\fg$ gives rise to a finite group $S(\fg)$ that coincides with the Weyl group of $\fg$ in case $\fg$ is a Lie algebra of classical type. In this paper, we compute the group $S(\fg)$ for Lie algebras of Cartan type and provide applications concerning weight space decompositions, the existence of generic tori and polynomial invariants. \end{abstract}

\maketitle

\section*{Introduction}
In the structure and representation theory of complex semi-simple Lie algebras, root systems and their Weyl groups play a fundamental r\^ole. For fields of positive characteristic the situation
is similar, provided one studies Lie algebras $\fg = \Lie(G)$ that are associated with a reductive group $G$. By contrast, the maximal tori of arbitrary restricted Lie algebras usually are no longer
conjugate, so there are many root systems available, whose structure and utility crucially depend on the initially chosen maximal torus. Accordingly, it is not clear which maximal tori are
appropriate for defining finite groups that take on the r\^ole of Weyl groups. In his article \cite{Pr2}, Premet studied this problem for the Jacobson-Witt algebras $W(n)$. He identified a
maximal torus $\ft_0 \subseteq W(n)$, whose automorphism group allowed to establish an analogue of the Chevalley restriction theorem on polynomial invariants.

One way to address the problem of choosing tori is the systematic study of schemes of tori. Roughly speaking, this geometric approach seeks to identify generic properties of tori and their
root systems. In the context of Lie algebras of complexity $\le 2$, that naturally arise in the study of restricted enveloping algebras of tame representation type, this method has produced
satisfactory results, cf.\ \cite{FV1,FV2}. A detailed investigation of the scheme of tori of maximal dimension $\mu(\fg)$ of a restricted Lie algebra $\fg$ in \cite{Fa2} led to the definition of a group $S(\fg)$ that turned out to coincide with the Weyl group in case $\fg$ is of classical type.

Aside from the analogues of the complex simple Lie algebras there are usually four additional classes of restricted simple Lie algebras, the so-called restricted Lie algebras of Cartan type.
These Lie algebras, which include the abovementioned Jacobson-Witt algebras, already differ significantly from their classical precursors in that their maximal tori are no longer conjugate. The
main object of the present paper is to identify the corresponding finite groups along with providing applications.

Our paper is organized as follows. After recalling basic results on the scheme $\cT_\fg$ of tori of a restricted Lie algebra $\fg$, we show in Section \ref{S:GTW} that the toral stabilizer
$S(\fg)$ of $\fg$ is isomorphic to the ``Weyl group" associated with any generic torus of $\fg$. For future reference, we also provide criteria for a given torus or Cartan subalgebra to be generic.
Section \ref{S:PTS} is concerned with the behavior of toral stabilizers under passage to subalgebras and factor algebras. In particular, we show that the toral stabilizer of a product of restricted
Lie algebras is the product of the toral stabilizers of the factors. Using basic results on tori and automorphism groups, due to Demu\v{s}kin \cite{De1,De2} and Wilson \cite{Wi1,Wi2},
respectively, we identify in Section \ref{S:AGT} the generic tori of the restricted Lie algebras of Cartan types $W,S$ and $H$. Accordingly, the Weyl groups of these tori coincide with the
toral stabilizers, whose computation relies on embeddings of certain subalgebras, particularly those of the form $W(\mu(\fg)) \hookrightarrow \fg$. Section \ref{S:TSCA} is devoted to the proof of our main result:

\bigskip

\begin{thm*} Suppose that $p\ge 3$ and let $\fg$ be a restricted Lie algebra of Cartan type. Then there is an isomorphism
\[ S(\fg) \cong \GL_{\mu(\fg)}(\FF_p)\]
and the scheme $\cT_\fg$ of embeddings of tori of maximal dimension is irreducible. \end{thm*}

\bigskip
\noindent
Along the way, a few other toral stabilizers need to be computed, notably those of the Poisson algebras $\cP(2r)$ with toral center.

The final section presents applications of our methods and results. With regard to weight space decompositions defined by tori of maximal dimension, the group $S(\fg)$ exhibits features analogous to those of Weyl groups: There is an action of $S(\fg)$ on the set of weights such that the dimensions of weight spaces are constant on the orbits. As a result, only three types of weight sets can occur for Lie algebras $\fg$ with $S(\fg) \cong \GL_{\mu(\fg)}(\FF_p)$. In particular, all root spaces of such Lie algebras have the same dimension.

In \cite[Chap.7]{St} Strade shows that the contact algebras $K(2r\!+\!1)$ possess infinitely many conjugacy classes of tori of dimension $\mu(K(2r\!+\!1))=r\!+\!1$. Using explicit information on Weyl groups and toral stabilizers of the aforementioned Poisson algebras, we prove that $K(2r\!+\!1)$ does not possess any generic tori.

Given a restricted Lie algebra $(\fg,[p])$, we denote by $S_\fg$ the variety of semi-simple elements. By modifying Premet's approach concerning $W(n)$, we finally establish the following analogue of the Chevalley restriction theorem: 

\bigskip

\begin{thm**} Suppose that $p\ge 3$ and let $(\fg,[p])$ be a restricted Lie algebra of Cartan type $W,S$ or $H$ with generic torus $\ft \subseteq \fg$. Then the restriction map induces an isomorphism $k[\bar{S}_\fg]^G \stackrel{\sim}{\lra} k[\ft]^{\GL_{\mu(\fg)}(\FF_p)}$. In particular, $k[\bar{S}_\fg]^G$ is a polynomial ring in $\mu(\fg)$ variables.
\end{thm**}

\bigskip

\begin{center}
{\bf Acknowledgements}
\end{center}

\bigskip We would like to thank the referee for providing additional references and useful comments.

\bigskip

\section{Generic Tori and Weyl Groups}\label{S:GTW}
Throughout, we shall be working over an algebraically closed field $k$ of characteristic $\ch(k) = p>0$. Unless mentioned otherwise, all algebras and modules are assumed to be finite-dimensional. Given a restricted Lie algebra $(\fg,[p])$, we let $G:= \Aut_p(\fg)^\circ$ be the identity component of its automorphism group. We denote by $\mu(\fg)$ and $\rk(\fg)$ the maximal dimension of all tori $\ft \subseteq \fg$ and the minimal dimension of all Cartan subalgebras $\fh \subseteq \fg$, respectively. According to \cite[(7.4)]{Fa2}, the set
\[\Tor(\fg) := \{ \ft \subseteq \fg \ ; \ \ft \ {\rm torus}, \ \dim_k\ft = \mu(\fg)\}\]
of tori of maximal dimension is locally closed within the Grassmanian ${\rm Gr}_{\mu(\fg)}(\fg)$ of $\fg$. Thanks to \cite[(1.6)]{FV1} and \cite[(3.5)]{Fa2}, the variety $\Tor(\fg)$ is
irreducible of dimension $\dim \Tor(\fg) = \dim_k\fg-\rk(\fg)$.

Let $M_k$ and ${\rm Ens}$ be the categories of commutative $k$-algebras (of arbitrary dimension) and sets, respectively. For a torus $\ft \in \Tor(\fg)$, we consider the
scheme $\cT_\fg : M_k \lra {\rm Ens}$, given by
\[ \cT_\fg (R) := \{ \varphi \in \Hom_p(\ft\!\otimes_k\!R,\fg\!\otimes_k\!R) \ ; \ \varphi \ \text{is a split injective $R$-linear map}\}\]
for every $R \in M_k$. Here $\Hom_p(\ft\!\otimes_k\!R,\fg\!\otimes_k\!R)$ denotes the set of homomorphisms of restricted $R$-Lie algebras and $\fg\!\otimes_k\!R$ carries the natural
structure of a restricted $R$-Lie algebra, with bracket and $p$-mapping defined via
\[ [x\otimes r,y\otimes s] := [x,y]\otimes rs \ \ \text{and} \ \ (x\otimes r)^{[p]} := x^{[p]}\otimes r^p \ \ \ \ \forall \ x,y \in \fg, \, r,s \in R,\]
respectively. According to \cite[(1.4),(1.6)]{FV1}, $\cT_\fg$ is a smooth affine scheme of dimension $\dim_k\fg-\rk(\fg)$.

In what follows, we shall be mainly concerned with the variety $\cT_\fg(k)$ of $k$-rational points of $\cT_\fg$, that is, the variety of embeddings $\varphi : \ft \lra \fg$ of restricted Lie
algebras. The automorphism groups $\Aut_p(\fg)$ and $\Aut_p(\ft)$ naturally act on $\cT_\fg(k)$ via
\[ g\dact \varphi := g\circ \varphi \ \ \text{and} \ \ h\dact \varphi := \varphi \circ h^{-1} \ \ \ \ \forall \ \varphi \in \cT_\fg(k), \, g \in \Aut_p(\fg), \, h \in \Aut_p(\ft), \]
respectively. Both actions commute, and $\Aut_p(\ft) \cong \GL_{\mu(\fg)}(\FF_p)$ is a finite group. We shall also consider the canonical action of $G$ on $\fg$ and write
\[ g.x := g(x) \ \ \ \ \forall \ g \in G, \, x \in \fg.\]
Let $\iota : \ft \hookrightarrow \fg$ be the standard embedding, defined by the inclusion $\ft \subseteq \fg$. Since $\cT_\fg(k)$ is smooth, its irreducible components coincide
with its connected components. Hence there exists exactly one irreducible component $\cX_\ft(k) \subseteq \cT_\fg(k)$ containing the inclusion $\iota$. We let
\[ S(\fg,\ft) := {\rm Stab}_{\Aut_p(\ft)}(\cX_\ft(k))\]
be the stabilizer of the component $\cX_\ft(k)$ in $\Aut_p(\ft)$. For future reference, we recall the following result, cf.\ \cite[(4.1)]{Fa2}:

\bigskip

\begin{Theorem} \label{GTWG1} Let $(\fg,[p])$ be a restricted Lie algebra, $\ft\subseteq \fg$ be a torus of dimension $\mu(\fg)$. Then the following statements hold:

{\rm (1)} \ $\cT_\fg(k)=\bigcup_{h\in \Aut_p(\ft)}h\dact\cX_\ft(k)$ and $\cT_\fg(k)$ has $[\Aut_p(\ft)\!:\!S(\fg,\ft)]$ irreducible components.

{\rm (2)} \ If $\ft' \subseteq \fg$ is another torus of dimension $\mu(\fg)$, then there exists an isomorphism $h : \ft \lra \ft'$ such that $S(\fg,\ft') = hS(\fg,\ft)h^{-1}$. \hfill $\square$
\end{Theorem}

\bigskip
\noindent
We let $\Nor_G(\ft)$ and $\Cent_G(\ft)$ be the {\it normalizer} and the {\it centralizer} of $\ft$ in $G$, respectively. In view of the above result, we make the following:

\bigskip

\begin{Definition} Let $(\fg,[p])$ be a restricted Lie algebra with automorphism group $G = \Aut_p(\fg)^\circ$, $\ft \subseteq \fg$ be a torus of dimension $\mu(\fg)$. Then
\[ S(\fg) := S(\fg,\ft)\]
is called the {\it toral stabilizer of $\fg$}. The group
\[ W(\fg,\ft):=\Nor_G(\ft)\slash \Cent_G(\ft)\]
is referred to as the {\it Weyl group of $\fg$ relative to $\ft$}. \end{Definition}

\bigskip
\noindent
The group $G$ naturally acts on the variety $\Tor(\fg)$ of tori of maximal dimension. Since $W(\fg,g\dact\ft) \cong W(\fg,\ft)$ for every $g \in G$, the Weyl group $W(\fg,\ft)$ only
depends on the orbit $G\dact\ft \subseteq \Tor(\fg)$.

\bigskip

\begin{Lemma} \label{GTWG2} Let $(\fg,[p])$ be a restricted Lie algebra with connected automorphism group $G := \Aut_p(\fg)^\circ$, $\ft \subseteq \fg$ be a torus of dimension
$\mu(\fg)$. Given $\varphi\in \cX_\ft(k)$, the map
\[\Theta_\varphi : \Nor_G(\varphi(\ft)) \lra \Aut_p(\ft) \ \ ; \ \ g \mapsto \varphi^{-1}\circ g \circ \varphi\]
is a homomorphism of groups that induces an injective homomorphism
\[ \bar{\Theta}_\varphi : W(\fg,\varphi(\ft)) \hookrightarrow S(\fg,\ft) \ \ ; \ \ \bar{g} \mapsto \varphi^{-1}\circ g \circ \varphi.\] \end{Lemma}

\begin{proof} Since the connected group $G$ acts morphically on $\cT_\fg(k)$, \cite[(8.2)]{Hu} ensures that $G$ stabilizes the connected component $\cX_\ft(k)$.  Let $g\in \Nor_G(\varphi(\ft))$. Setting $h:=\Theta_\varphi(g)$, we consider the automorphism
\[ \rho_h : \cT_\fg(k) \lra \cT_\fg(k) \ \ ; \ \ \psi \mapsto \psi\circ h\]
of the variety $\cT_\fg(k)$. Then $\rho_h(\cX_\ft(k))$ is a connected component of $\cT_\fg(k)$. Since
\[ \rho_h(\varphi) = \varphi \circ \Theta_\varphi(g) = g \circ \varphi \in \cX_\ft(k),\]
we obtain $\rho_h(\cX_\ft(k))=\cX_\ft(k)$. As a result, $\Theta_\varphi(g)\in S(\fg,\ft)$. Directly from the definition we obtain the identity $\ker \Theta_\varphi=\Cent_G(\varphi(\ft))$,
so that the induced map $\bar{\Theta}_\varphi : W(\fg,\varphi(\ft)) \lra S(\fg,\ft)$ is actually injective. \end{proof}

\bigskip
\noindent
By virtue of Theorem \ref{GTWG1}, every torus $\ft'\subseteq \fg$ of dimension $\mu(\fg)$ is of the form $\ft'=\varphi(\ft)$ for some $\varphi\in \cX_\ft(k)$. The foregoing Lemma then
provides an embedding
\[\bar{\Theta}_\varphi : W(\fg,\ft')\hookrightarrow S(\fg,\ft),\]
sending the Weyl group $W(\fg,\ft')$ of $\ft'$ into the toral stabilizer $S(\fg,\ft)$. If $\fg$ is of classical type, then this map is an isomorphism, see \cite[(4.6)]{Fa2} or Theorem \ref{GTWG5} below.

Our main goal of this section is to provide a criterion ensuring that $W(\fg,\ft) \cong S(\fg)$ for some torus $\ft \in \Tor(\fg)$.

\bigskip

\begin{Definition} A torus $\ft \in \Tor(\fg)$ is called {\it generic} if the orbit $G\dact\ft$ is a dense subset of $\Tor(\fg)$. \end{Definition}

\bigskip
\noindent
Let $\fg := \Lie(G)$ be the Lie algebra of a connected algebraic group. In view of \cite[(13.3),(13.5)]{Hu0} (see also \cite[(4.4)]{Fa2}), any two tori of $\fg$ of maximal dimension are conjugate via the adjoint action. Accordingly, each of these tori is generic.

\bigskip

\begin{Lemma} \label{GTWG3} Let $(\fg,[p])$ be a restricted Lie algebra with connected automorphism group $G:=\Aut_p(G)^\circ$, $\ft \subseteq \fg$ be a torus of dimension $\mu(\fg)$.
Then
\[\dim G.\varphi(\ft)=\dim G\dact \varphi+\mu(\fg)\]
 for every $\varphi\in \cX_\ft(k)$. \end{Lemma}

\begin{proof} We consider the morphism
\[ \omega : \left \{ \begin{array}{ccc} \cX_\ft(k) \times \ft  & \lra & \fg \\
                                                                     (\varphi,t) & \mapsto & \varphi(t)
\end{array} \right.\]
of irreducible varieties. Given $t \in \ft$, we let $(kt)_p$ denote the $p$-subalgebra of $\fg$ generated by $kt$. Lemma 3.2 of \cite{Fa2} provides a non-empty open subset $U\subseteq \ft$
with $\ft=(kt)_p$ for every $t\in U$.  Let $(\varphi,t)$ be an element of $\cX_\ft(k)\!\times\!U$. If $(\psi,s)\in \omega^{-1}(\omega(\varphi,t))$, then $\psi(s)=\varphi(t)$, so that
\[\varphi(t)\in \im \psi.\]
By choice of $t$, this implies that $\im \varphi \subseteq \im\psi$, whence $\im \varphi=\im\psi$. Hence there exists an element $h\in \Aut_p(\ft)$ with
\[\psi=\varphi\circ h^{-1}.\]
In particular, $\varphi(h^{-1}(s))=\psi(s)=\varphi(t)$, so that $s=h(t)$. As a result, the fiber $\omega^{-1}(\omega(\varphi,t))\subseteq \{(\varphi\circ h^{-1},h(t)) \ ; \ h\in
\Aut_p(\ft)\}$ is finite.

Now let $\varphi\in \cX_\ft(k)$ be arbitrary. Then $\omega$ induces a surjective morphism
\[\omega:(G\dact\varphi)\times \ft \lra G.\varphi(\ft)\]
of irreducible varieties. By the above argument, the fibers of this morphism are finite on the dense open subset $(G\dact\varphi)\!\times\! U$ of $(G\dact\varphi)\!\times\!\ft$. By virtue of
\cite[I.\S8, Corollary 1]{Mu}, we obtain
\begin{align*}\dim G.\varphi(\ft) &=\dim ((G\dact\varphi)\!\times\!\ft) = \dim G\dact\varphi +\dim_k\ft\\
&=\dim G\dact\varphi +\mu(\fg),
\end{align*}
as desired. \end{proof}

\bigskip
\noindent
Let $S_\fg$ denote the set of semi-simple elements of $\fg$. The following result relates dense orbits in $\Tor(\fg)$ and $\cX_\ft(k)$ to dense $G$-saturations of tori in $\overline{S}_\fg$
relative to the canonical action of $G$ on $\fg$.

\bigskip

\begin{Proposition} \label{GTWG4} Let $(\fg,[p])$ be a restricted Lie algebra, $G = \Aut_p(\fg)^\circ$ be its connected automorphism group, and $\ft \subseteq \fg$ be a torus of maximal
dimension. For a torus  $\ft_0 \in \Tor(\fg)$, the following statements are equivalent:

{\rm (1)} \ The torus $\ft_0$ is generic.

{\rm (2)} \ There exists $\varphi_0 \in \cX_\ft(k)$ such that $\varphi_0(\ft) = \ft_0$ and $\overline{G\dact\varphi_0} = \cX_\ft(k)$.

{\rm (3)} \ $G.\ft_0$ is a dense subset of $\overline{S}_\fg$. \end{Proposition}

\begin{proof} (1) $\Rightarrow$ (2). According to \cite[p.~4226]{Fa2}, the map
  \[ \zeta_\ft : \cX_\ft(k) \lra \Tor(\fg) \ \ ; \ \ \varphi \mapsto \varphi(\ft)\]
is a surjective morphism with finite fibers between two irreducible varieties of dimension $\dim_k\fg -\rk(\fg)$. Let $\varphi \in \zeta_{\ft}^{-1}(\ft_0)$, so that $\varphi_0(\ft) = \ft_0$. Since
$\zeta_\ft$ is $G$-equivariant, our present assumption implies that the morphism $\zeta_\ft|_{\overline{G\dact\varphi_0}} :  \overline{G\dact\varphi_0} \lra \Tor(\fg)$ is dominant.
Consequently,
\[ \dim \overline{G\dact\varphi_0} = \dim \Tor(\fg) = \dim \cX_\ft(k),\]
so that $\cX_\ft(k) = \overline{G\dact\varphi_0}$.

(2) $\Rightarrow$ (3). According to \cite[(4.1)]{Fa2} and \cite[(3.7)]{Fa2} (see also \cite[Theorem~2(iii)]{Pr1}), the variety $\bar{S}_\fg$ is irreducible, and of dimension
\[ \dim \bar{S}_\fg = \dim \cX_\ft(k) + \mu(\fg).\]
Lemma \ref{GTWG3} now implies
\[ \dim G.\ft_0 = \dim G\dact\varphi_0 + \mu(\fg) = \dim \overline{S}_\fg,\]
so that $G.\ft_0$ is a dense subset of the irreducible variety $\overline{S}_\fg$.

(3) $\Rightarrow$ (1). Thanks to Theorem \ref{GTWG1}, there exists an element $\varphi_0 \in \cX_\ft(k)$ with $\ft_0 = \varphi_0(\ft)$. According to Lemma \ref{GTWG3}, we have
\[ \dim G\dact\varphi_0 = \dim G.\ft_0 - \mu(\fg) = \dim \bar{S}_\fg - \mu(\fg) = \dim \cX_\ft(k),\]
so that the orbit $G\dact\varphi_0$ lies dense in $\cX_\ft(k)$. Hence
\[ \Tor(\fg) = \zeta_\ft(\cX_\ft(k)) = \zeta_\ft(\overline{G\dact\varphi_0}) \subseteq \overline{\zeta_\ft(G\dact\varphi_0)} = \overline{G\dact\ft_0},\]
and $\ft_0$ is a generic torus. \end{proof}

\bigskip

\begin{Theorem} \label{GTWG5} Let $(\fg,[p])$ be a restricted Lie algebra with connected automorphism group $G := \Aut_p(\fg)^\circ$. If $\ft \in \Tor(\fg)$ is a generic torus,
 then the map
 \[ \Nor_G(\ft) \lra \Aut_p(\ft) \ \ ; \ \ g \mapsto g|_\ft \]
 induces an isomorphism $ W(\fg,\ft) \cong S(\fg,\ft) $ of groups. \end{Theorem}

\begin{proof} Thanks to Proposition \ref{GTWG4}, there exists $\varphi_0 \in \cX_\ft(k)$ with $\varphi_0(\ft) = \ft$ such that the orbit $G\dact\varphi_0$ lies dense in $\cX_\ft(k)$. It follows
that $G\dact\varphi_0$ is an open subset of the irreducible variety $\cX_\ft(k)$. Consequently, $G\dact\varphi_0$ is the unique dense orbit of $\cX_\ft(k)$.

Given $h\in S(\fg,\ft)$, we consider the automorphism
\[ w_h : \cX_\ft(k) \lra \cX_\ft(k) \ \ ; \ \ \varphi\mapsto \varphi\circ h\]
of the variety $\cX_\ft(k)$. Then $w_h$ is $G$-equivariant, so that
\[ w_h(G\dact\varphi_0) = G\dact w_h(\varphi_0)\]
is a dense orbit of $\cX_\ft(k)$. By our observation above, we have $G\dact w_h(\varphi_0) = G\dact\varphi_0$, and there exists an element $g\in G$ with
\begin{equation} \label{astre}
w_h(\varphi_0)=g\circ \varphi_0.
\end{equation}
In particular, we obtain
\[g(\varphi_0(t)) = \varphi_0(h(t)) \in \varphi_0(\ft)=\ft\]
for every $t\in \ft$, whence $g\in \Nor_G(\varphi_0(\ft))=\Nor_G(\ft)$. Consequently, (\ref{astre}) implies $h = \Theta_{\varphi_0}(g)$, so that $\Theta_{\varphi_0}$ is surjective. Thanks to Lemma
\ref{GTWG2}, the homomorphism $\bar{\Theta}_{\varphi_0} : W(\fg,\ft) \lra S(\fg,\ft)$ is bijective, so that the groups $W(\fg,\ft)$ and $S(\fg,\ft)$ have the same order. We may now apply Lemma \ref{GTWG2} to the canonical embedding $\iota : \ft \hookrightarrow \fg$ to see that $\bar{\Theta}_{\iota}$ is the desired
isomorphism.\end{proof}

\bigskip
\noindent
We provide another density result which compensates for the lack of a generic torus to some extent. Given a torus $\ft \subseteq \fg$, we let $C_\fg(\ft)$ be the {\it centralizer} of $\ft$ in $\fg$.

\bigskip

\begin{Lemma} \label{GTWG8} Let $(\fg,[p])$ be a restricted Lie algebra. Then $ \fT := \bigcup_{\ft \in \Tor(\fg)} \ft$ is a dense subset of $\bar{S}_\fg$. \end{Lemma}

\begin{proof}
We put $A = k[\cX_\ft]$ and $\tilde{\fg} = \fg \otimes_k A$. Let $j:\ft \hookrightarrow \tilde{\fg}$ be the universal embedding. The torus $\ft$ acts on $\tilde{\fg}$ via $j$; we denote by $\tilde{\fg}_0 \subseteq \tilde{\fg}$ the corresponding zero weight space. For any homomorphism $x \in \Spec_k(A)(k)$, let $\varphi_x := (\id_{\fg} \otimes x) \circ j \in \cT_\fg(k)$ denote the corresponding embedding. Then $\fg_0(x) = (\id_\fg \otimes x) (\tilde{\fg}_0)$ is the zero weight space relative to the torus $\varphi_x(\ft)$. Consider the variety
\[\cV = \{ (x,y) \in \cX_{\ft}(k) \times \fg\ ; \ \exists\, v \in \tilde{\fg}_0\ \mbox{with}\ (\id_{\fg} \otimes x) (v) = y \}. \]
Owing to \cite[(3.4)]{Fa2}, the projection $\pi_{\fg}: \cV \lra \fg$ is a dominant morphism with image
\[\im \pi_{\fg} = \bigcup_{x \in \cX_\ft(k)} \tilde{\fg}_0(x) = \bigcup_{x \in \cX_{\ft}(k)} C_{\fg}(\varphi_x(\ft)).\]
Also, for $n = \dim_k \fg$, the morphism $\zeta: \fg \lra \bar{S}_\fg \ ; \ y \mapsto y^{[p]^n}$ is dominant. Thus, $\zeta \circ \pi_{\fg}$ is dominant, and we have
\begin{equation} \label{imzetarondpi}
\im \zeta \circ \pi_{\fg} = \bigcup_{x \in \cX_\ft(k)} C_{\fg}(\varphi_x(\ft))^{[p]^n} = \bigcup_{x \in \cX_\ft(k)} \varphi_x(\ft).
\end{equation}
Since the right hand set is dense in $\bar{S}_\fg$ and contained in $\fT$, our assertion follows.\end{proof}

\bigskip
\begin{Remark}
In the proof, we used the fact that the morphism $\fg \lra \fg \ ; \ y \mapsto y^{[p]^n}$ is dominant for $n = \dim_k \fg$. In fact, there are lower values of $n$ for which this property holds, see the Remark p. \pageref{seven} and Corollary \ref{IL6}.
\end{Remark}

\bigskip
\noindent
We next record analogous properties of Cartan subalgebras. Thanks to \cite[(7.2),(7.4)]{Fa2}
\[ \Car(\fg) := \{ \fh \subseteq \fg \ ; \ \fh \ \text{Cartan subalgebra}, \ \dim_k \fh = \rk(\fg)\}\]
is a locally closed subset of the Grassmann variety ${\rm Gr}_{\rk(\fg)}(\fg)$. In view of \cite[(8.2)]{Fa2}, the variety $\Car(\fg)$ is irreducible of dimension $\dim_k\fg\!-\!\rk(\fg)$. Note that $G$ naturally acts on $\Car(\fg)$.

By general theory, the Cartan subalgebras of $\fg$ are precisely the centralizers of the maximal tori, see for instance \cite[(II.4.1)]{SF}. Moreover, $C_\fg(\ft) \in \Car(\fg)$ for every $\ft \in \Tor(\fg)$, cf.\ \cite[(3.5(1))]{Fa2}.

\bigskip
\begin{Remark} \label{seven}
Let $\ell:=\rk(\fg)-\mu(\fg)$, and $\zeta_\ell: \fg \lra \fg \ ; \ y \mapsto y^{[p]^\ell}$. Then $\zeta_\ell$ induces a dominant morphism $\fg \to \bar{S}_\fg$. Indeed, since $\pi_{\fg}: \cV \lra \fg$ is dominant, it is enough to show that $\zeta_\ell \circ \pi_{\fg}$ is dominant. We first check it under for $\fg$ nilpotent. Then the set $N(\fg)$ of $p$-unipotent elements form a $p$-subalgebra of $\fg$ of dimension $\ell$. Hence, for all $x_n \in N(\fg)$, we have $x_n^{[p]^\ell} = 0$. Furthermore, $\fg$ affords a unique torus $\ft$ of dimension $\mu(\fg)$. Then it is then easy to see that $x^{[p]^\ell} \in \ft$, whence $\fg^{[p]^\ell} = \ft = S_\fg$ as we needed. Now when $\fg$ is arbitrary, we have
\[\im \zeta_\ell \circ \pi_{\fg} = \bigcup_{x \in \cX_\ft(k)} C_{\fg}(\varphi_x(\ft))^{[p]^\ell}.\]
Since each $\fh:=C_{\fg}(\varphi_x(\ft))$ is a nilpotent subalgebra of dimension $\rk(\fg)$ with $\mu(\fh) = \mu(\fg)$, we have $\fh^{[p]^\ell} = \varphi_x(\ft)$. Comparing with Equation (\ref{imzetarondpi}) yields $\im \zeta_\ell \circ \pi_{\fg} = \im \zeta \circ \pi_{\fg}$, a dense subset of $\bar{S}_\fg$.
\end{Remark}

\bigskip

\begin{Definition} A Cartan subalgebra $\fh \in \Car(\fg)$ is called {\it generic} if the orbit $G\dact \fh$ is dense in $\Car(\fg)$. \end{Definition}

\bigskip
\noindent
As before, the Cartan subalgebras of the Lie algebras that are associated with algebraic groups are generic.

\bigskip

\begin{Proposition} \label{GTWG6} Let $(\fg,[p])$ be a restricted Lie algebra with connected automorphism group $G := \Aut_p(\fg)^\circ$. Suppose that $\ft_0 \in \Tor(\fg)$ is a torus, and let
$\fh_0 := C_\fg(\ft_0) \in \Car(\fg)$ be the corresponding Cartan subalgebra. Then the following statements are equivalent:

{\rm (1)} \ The Cartan subalgebra $\fh_0$ is generic.

{\rm (2)} \ The torus $\ft_0$ is generic.

{\rm (3)} \ The $G$-saturation $G.\fh_0 \subseteq \fg$ is a dense subset of $\fg$.\end{Proposition}

\begin{proof} (1) $\Rightarrow$ (2). A consecutive application of \cite[(9.3)]{Fa2} and \cite[(8.2)]{Fa2} shows that
\[ C_\fg : \Tor(\fg) \lra \Car(\fg) \ \ ; \ \ \ft \mapsto C_\fg(\ft)\]
is an injective morphism, whose image $\im C_\fg = \{\fh \in \Car(\fg) \ ; \ \mu(\fh) = \mu(\fg)\}$ is a dense open subset of $\Car(\fg)$. Since $C_\fg$ is $G$-equivariant, we obtain
\[ C_\fg(G\dact\ft_0) = G\dact\fh_0,\]
so that the morphism $C_\fg : \overline{G\dact\ft_0} \lra \Car(\fg)$ is dominant. Consequently,
\[ \dim \overline{G\dact\ft_0} = \dim \Car(\fg) = \dim \Tor(\fg).\]
As $\Tor(\fg)$ is irreducible, we conclude $\overline{G\dact\ft_0} = \Tor(\fg)$, showing that $\ft_0$ is a generic torus.

(2) $\Rightarrow$ (3). It follows from \cite[(1.6)]{Fa2} that the variety
\[ \cV := \{(x,\varphi) \in \fg\times\cX_{\ft_0}(k) \ ; \ x \in C_\fg(\varphi(\ft_0))\}\]
is irreducible. Moreover, \cite[(3.1)]{Fa2} implies that the morphism
\[ \pi_\fg : \cV \lra \fg \ \ ; \ \ (x,\varphi) \mapsto x\]
is dominant. Since $\ft_0$ is generic, Proposition \ref{GTWG4} provides an element $\varphi_0 \in \cX_{\ft_0}(k)$ with image $\varphi_0(\ft_0) = \ft_0$ and such that
$\overline{G\dact\varphi_0} = \cX_{\ft_0}(k)$. Hence $G\dact\varphi_0$ is open in $\cX_{\ft_0}(k)$, and
\[ \cO := (\fg\times G\dact\varphi_0)\cap\cV\]
is open in $\cV$. Thus, $\cO$ is a dense subset of $\cV$, and we obtain
\[ \fg = \overline{\pi_\fg(\cV)} =  \overline{\pi_\fg(\overline{\cO})} \subseteq \overline{\pi_\fg(\cO)} = \overline{G.C_\fg(\varphi_0(\ft_0))} = \overline{G.\fh_0},\]
as desired.

(3) $\Rightarrow$ (1). Suppose that $\fg = \overline{G.\fh_0}$. Setting $\ell := \dim_k\fg$, we recall that
\[ f : \fg \lra \fg \ \ ; \ \ x \mapsto x^{[p]^\ell}\]
is a $G$-equivariant morphism, whose image coincides with $S_\fg$. Moreover,
\[ f(\fh_0) = \fh_0\cap S_\fg =\ft_0,\]
so that $f(G.\fh_0) = G.\ft_0$. Consequently,
\[ \overline{S_\fg} = \overline{f(\fg)} = \overline{f(\overline{G.\fh_0})} \subseteq \overline{f(G.\fh_0)} = \overline{G.\ft_0},\]
and Proposition \ref{GTWG4} implies that $\ft_0$ is generic. Thus,
\[ \Car(\fg) = \overline{C_\fg(\Tor(\fg))} \subseteq \overline{C_\fg(G\dact\ft_0)} = \overline{G\dact\fh_0},\]
proving that $\fh_0$ is a generic Cartan subalgebra of $\fg$.  \end{proof}

\bigskip

\begin{Lemma} \label{GTWG7} Let $(\fg,[p])$ be a restricted Lie algebra, $V \subseteq \fg$ be a subspace of $\fg$.

{\rm (1)} \ $\Tor(\fg)^V := \{ \ft \in \Tor(\fg) \ ; \ \ft \cap V = (0)\}$ is open in $\Tor(\fg)$.

{\rm (2)} \ $\Tor(\fg)_V := \{ \ft \in \Tor(\fg) \ ; \ \ft \subseteq V\}$ is closed in $\Tor(\fg)$.

{\rm (3)} \ $\Car(\fg)^V := \{ \fh \in \Car(\fg) \ ; \ \fh \cap V = (0)\}$ is open in $\Car(\fg)$.

{\rm (4)} \ $\Car(\fg)_V := \{ \fh \in \Car(\fg) \ ; \ \fh \subseteq V\}$ is closed in $\Car(\fg)$. \end{Lemma}

\begin{proof} Let $d \in \NN$. According to \cite[(7.3)]{Fa2}, the map
\[ f_{d,V} : {\rm Gr}_d(\fg) \lra \NN_0 \ \ ; \ \ X \mapsto \dim_kX\cap V\]
is upper semicontinuous. Consequently, $\Tor(\fg)^V = \Tor(\fg)\setminus f^{-1}_{\mu(\fg),V}(\{n \in \NN_0 \ ; \ n \ge 1\})$ is open in $\Tor(\fg)$. By the same token, $\Tor(\fg)_V =
\Tor(\fg)\cap f_{\mu(\fg),V}^{-1}(\{n \in \NN_0 \ ; \ n \ge \mu(\fg)\})$ is closed. The proofs of (3) and (4) follow analogously. \end{proof}

\bigskip

\section{Basic Properties of Toral Stabilizers}\label{S:PTS}
Our computations of toral stabilizers necessitate information concerning their behavior with respect to inclusions and passage to factor algebras. Let $(\fg,[p])$ be a restricted Lie algebra with
torus $\ft \subseteq \fg$ of maximal dimension. Recall that $S(\fg,\ft)$ is a subgroup of $\Aut_p(\ft)$.

\bigskip

\begin{Lemma} \label{BPTS1} Let $\fh \subseteq \fg$ be a $p$-subalgebra of the restricted Lie algebra $(\fg,[p])$.

 {\rm (1)} \ If $\mu(\fh) = \mu(\fg)$, then $S(\fh,\ft) \subseteq S(\fg,\ft)$ for every torus $\ft \subseteq \fh$ of dimension $\mu(\fh)$.

 {\rm (2)} \ If $\fh \unlhd \fg$ is an ideal such that $\fg/\fh$ is $p$-unipotent, then $\mu(\fh)=\mu(\fg)$ and $S(\fh,\ft) = S(\fg,\ft)$ for every torus $\ft \subseteq \fh$ of dimension
 $\mu(\fh)$.\end{Lemma}

\begin{proof} (1) Let $\ft \subseteq \fh$ be a torus of dimension $\mu(\fh)$. The canonical inclusion $\iota : \fh \hookrightarrow \fg$ induces an injective morphism
\[ \iota_\ast : \cT_\fh(k) \lra \cT_\fg(k) \ \ ; \ \ \varphi \mapsto \iota\circ \varphi.\]
Consequently, $\iota_\ast(\cX_\ft(k)) \subseteq \cT_\fg(k)$ is a connected subspace containing the canonical inclusion $\iota|_\ft : \ft \hookrightarrow \fg$, so that
\[ \iota_\ast(\cX_\ft(k)) \subseteq \cX'_\ft(k),\]
where the latter space denotes the connected component of $\cT_\fg(k)$ containing $\iota|_\ft$. Let $h$ be an element of $S(\fh,\ft)$. Denoting the inclusion $\ft \hookrightarrow \fh$ by
$\tilde{\iota}$, we obtain
\[ h\dact \iota|_\ft = \iota_\ast(h\dact \tilde{\iota}) \in \cX'_\ft(k).\]
Since $\varphi \mapsto h\dact \varphi$ is a morphism, we see that $h\dact\cX'_\ft(k) \subseteq \cX'_\ft(k)$. This implies $h \in S(\fg,\ft)$, as desired.

(2) Since $\fg/\fh$ is $p$-unipotent, we have $\mu(\fg/\fh)=0$ and \cite[(3.3)]{FV2} implies $\mu(\fh) = \mu(\fg)$. Let $\pi : \fg \lra \fg/\fh$ be the canonical projection. Given
$\varphi \in \cT_\fg(k)$, the $p$-subalgebra $(\pi\circ\varphi)(\ft)\subseteq \fg/\fh$ is a torus, so that $(\pi\circ \varphi)(\ft)=(0)$ and $\varphi(\ft)\subseteq \ker\pi = \fh$. Consequently, $\varphi$ factors through $\fh$, and the morphism
\[ \iota_\ast : \cT_\fh(k) \lra \cT_\fg(k)\]
is bijective. As a result, $\iota_\ast(\cX_\ft(k))=\cX'_\ft(k)$ and $S(\fg,\ft) = S(\fh,\ft)$. \end{proof}

\bigskip

\begin{Corollary} \label{BPTS2} Let $\ft \subseteq \fg$ be a torus of dimension $\mu(\fg)$, $\fs \subseteq \ft$ be a subtorus. Then we have a natural inclusion
\[ S(C_\fg(\fs),\ft) \subseteq \{ g \in S(\fg,\ft) \ ; \ g|_\fs = \id_\fs\}.\]
\end{Corollary}

\begin{proof} We consider the $p$-subalgebra $\fh := C_\fg(\fs)$. By definition, $\mu(\fh) = \mu(\fg)$, and Lemma \ref{BPTS1} provides an inclusion $S(\fh,\ft) \subseteq S(\fg,\ft)$.
In view of \cite[(5.11)]{Fa2}, each element of $S(\fh,\ft)$ belongs to the centralizer of $\fs$ in $\Aut_p(\ft)$. \end{proof}

\bigskip
\noindent
For a direct sum $\ft_1 \oplus \ft_2$ of vector spaces, we denote by $\iota_{\ft_i}$ and $\pi_{\ft_i}$ the canonical injection from $\ft_i$ and the projection onto $\ft_i$, respectively. Recall
the identification
\begin{equation} \label{IdentifyAutGroups}
\renewcommand{\arraystretch}{1.7}
\End_k(\ft_1 \oplus \ft_2) \equiv \left( \begin{array}{c|c} \End_k(\ft_1) & \Hom_k(\ft_2,\ft_1) \\ \hline
\Hom_k(\ft_1,\ft_2) & \End_k(\ft_2) \end{array} \right),
\end{equation}
where the maps from $\End_k(\ft_1 \oplus \ft_2)$ to any of the blocks arise via composition with $\iota$'s and $\pi$'s: For example, the projection $\End_k(\ft_1 \oplus \ft_2)
\twoheadrightarrow \Hom_k(\ft_2 , \ft_1)$ is given by $g \mapsto  \pi_{\ft_1} \circ g \circ \iota_{\ft_2}$.

The groups $\GL(\ft_1)$ and $\GL(\ft_2)$ canonically act on $\Hom_k(\ft_2,\ft_1)$ from the left and right, respectively. The semi-direct product $\GL(\ft_2) \ltimes \Hom_k(\ft_2,\ft_1)$
can be identified with a subgroup of $\GL(\ft_1 \oplus \ft_2)$, namely:
\begin{equation*}
\renewcommand{\arraystretch}{1.7}
\GL(\ft_2) \ltimes \Hom_k(\ft_2,\ft_1) \equiv \left( \begin{array}{c|c} I_{\ft_1} & \Hom_k(\ft_2,\ft_1) \\ \hline 0 & \GL(\ft_2) \end{array} \right).
\end{equation*}
(Note that the sizes of boxes are given by the dimensions of the relevant subspaces of $\ft$, so that there is no entry in the first row and column for $t_1 = (0)$, say.)

\bigskip
\noindent
For tori $\ft,\ft'$, let $\Lie_p(\ft,\ft')$ denote the additive group of restricted homomorphisms $\ft \to \ft'$. It affords a natural right action by $\Aut_p(\ft)$ and a natural left action by $\Aut_p(\ft')$. We first generalize a basic result from \cite{Fa2}.

\bigskip

\begin{Proposition} \label{BPTS3} Let $(\fg,[p])$ be a restricted Lie algebra, $\fn \unlhd \fg$ be a $p$-ideal. Then the following statements hold:

 {\rm (1)} \ If $\ft \subseteq \fg$ is a torus of maximal dimension, then $\ft\cap\fn$ and $(\ft+\fn)/\fn$ are tori of maximal dimension of $\fn$ and $\fg/\fn$, respectively.

 {\rm (2)} \ The restriction map $\Hom_k(\ft,\fg) \lra \Hom_k(\ft\cap \fn,\fg)$ induces a morphism ${\rm res} : \cX_\ft(k) \lra \cX_{\ft\cap\fn}(k)$.

 {\rm (3)} \ Let $\ft' \subseteq \ft$ be a torus such that $\ft = (\ft\cap\fn)\oplus\ft'$. Then the canonical identification (\ref{IdentifyAutGroups}) induces an injective homomorphism
\[ \renewcommand{\arraystretch}{1.7} S(\fg , \ft) \hookrightarrow \left(\begin{array}{c|c} S(\fn,\ft\cap\fn) & \Lie_p(\ft',\ft\cap\fn)\\ \hline 0 & S(\fg/\fn,\ft')\end{array} \right)\]
of finite groups. \end{Proposition}

\begin{proof} (1) According to \cite[(3.3)]{FV2}, the $p$-subalgebra $(\ft+\fn)/\fn \subseteq \fg/\fn$ is a torus of $\fg/\fn$ of maximal dimension. By the same token, we have
\[ \mu(\fg) = \mu(\fn)+\mu(\fg/\fn).\]
Since $\dim_k\ft = \dim_k(\ft\cap\fn)+\dim_k(\ft+\fn)/\fn$, it follows that $\ft\cap\fn$ is a torus of $\fn$ of maximal dimension.

(2) Let $\iota : \ft \hookrightarrow \fg$ be the canonical embedding. We first verify

\medskip
($\ast$) \ {\it We have $\varphi(\ft\cap\fn) \subseteq \fn$ for every $\varphi \in \cX_\ft(k)$ and $g(\ft\cap\fn) = \ft\cap\fn$ for every $g \in S(\fg,\ft)$}.

\smallskip
\noindent
According to \cite[(5.2)]{Fa2}, there exists a subtorus $\ft_\fn \subseteq \ft$ such that $\varphi^{-1}(\fn) = \ft_\fn$ for every $\varphi \in \cX_\ft(k)$. Specializing $\varphi = \iota$, we
obtain $\ft_\fn = \ft\cap\fn$, whence $\varphi(\ft\cap\fn) \subseteq \fn$ for every $\varphi \in \cX_\ft(k)$.

Now let $g \in S(\fg,\ft)$. Then we have $g = g^{-1}\dact\iota \in \cX_\ft(k)$, so that $g(\ft\cap\fn) \subseteq \ft\cap \fn$. \hfill $\diamond$

\medskip
\noindent
Owing to ($\ast$), the canonical restriction map induces a morphism
\[ {\rm res} : \cX_\ft(k) \lra \cT_{\fn}(k) \ \ ; \ \ \varphi \mapsto \varphi|_{\ft\cap \fn}\]
of affine varieties. As ${\rm res}$ is continuous, we have $\im {\rm res} \subseteq \cX_{\ft\cap \fn}(k)$.

(3) \ Let $g$ be an element of $S(\fg,\ft)$. Owing to ($\ast$), we have $g(\ft\cap\fn) = \ft\cap \fn$. Thus,
\[  \iota|_{\ft\cap\fn} \circ g|_{\ft\cap\fn} = (\iota \circ g)|_{\ft\cap \fn} = {\rm res}(\iota\circ g) \in \cX_{\ft\cap\fn}(k),\]
proving $g|_{\ft\cap\fn} \in S(\fn,\ft\cap \fn)$.

According to \cite[(5.3)]{Fa2}, the natural projection $\pi : \fg \lra \fg/\fn$ induces a morphism
\[ \Pi : \cX_\ft(k) \lra \cX_{\ft'}(k) \ \ ; \ \ \varphi \mapsto \pi\circ \varphi|_{\ft'}\]
of affine varieties. Property ($\ast$) implies that each $g \in S(\fg,\ft)$ induces a unique element $\bar{g} \in \Aut_p((\ft+\fn)/\fn)$ such that $\bar{g}\circ \pi|_{\ft} = \pi \circ g$. Since
$\pi$ restricts to an isomorphism $\ft' \lra (\ft+\fn)/\fn$, there is a unique automorphism $\tilde{g} \in \Aut_p(\ft')$  such that $\bar{g}\circ \pi|_{\ft'} = \pi|_{\ft'}\circ \tilde{g}$.
Note that $\pi|_{\ft'} : \ft' \hookrightarrow \fg/\fn$ is the canonical inclusion, defining the connected component $\cX_{\ft'}(k) \subseteq \cT_{\fg/\fn}(k)$. Hence we obtain for $g \in
S(\fg,\ft)$
\[ \tilde{g}\dact \pi|_{\ft'} = \pi|_{\ft'}\circ \tilde{g}^{-1} = \pi \circ g^{-1}|_{\ft'} = \Pi(g^{-1}) \in \cX_{\ft'}(k),\]
so that $\tilde{g} \in S(\fg/\fn,\ft')$. Since $g \mapsto \tilde{g}$ defines a homomorphism of groups, we obtain the desired embedding of finite groups. \end{proof}

\bigskip

\begin{Remarks} (1) The first part of the Proposition is well-known, see for example \cite[Theorem 4.3]{Win} for maximal tori and \cite[Lemma 1.7.2]{BW} for tori of maximal dimension.

(2) Let $\ft' \subseteq \ft$ be tori. Let $\{t_1,\ldots, t_n\}$ be a toral basis of $\ft$ and denote by $V$ the $\FF_p$-subspace of $\ft$ generated by $\{t_1,\ldots, t_n\}$.
Since $V$ is the set of toral elements of $\ft$, there exists a subspace $V' \subseteq V$, whose $k$-span coincides with $\ft'$. Choose a complement $W'$ of $V'$ in $V$. Then its $k$-span $\fs$ is a subtorus of $\ft$ such that $\ft = \ft' \oplus \fs$.

(3) One can show that the restriction morphism ${\rm res} : \cX_\ft(k) \lra \cX_{\ft\cap\fn}(k)$ is smooth. In particular, ${\rm res}$ is open and dominant. \end{Remarks}

\bigskip

\begin{Corollary} \label{BPTS4} Let $\fg_1,\fg_2$ be restricted Lie algebras and $\ft_i \subseteq \fg_i$ be tori of maximal dimension for $i \in \{ 1 , 2 \}$. Then $\ft_1 \oplus \ft_2 \subseteq \fg_1 \oplus \fg_2$ is a torus of maximal dimension, and
\[S(\fg_1 \oplus \fg_2 , \ft_1 \oplus \ft_2) \simeq S(\fg_1 , \ft_1)\!\times\! S(\fg_2 , \ft_2).\] \end{Corollary}

\begin{proof} We put $\fg := \fg_1 \oplus \fg_2$ and $\ft := \ft_1 \oplus \ft_2$, and denote by $\iota : \ft \hookrightarrow \fg$ the canonical embedding. Let $\cT_\fg$ be the scheme of
embeddings of $\fg$ and $\cX \subseteq \cT_\fg$ be the connected component containing $\iota$. For $i \in \{ 1,2 \}$, we denote similarly $\iota_i : \ft_i \hookrightarrow \fg_i$ and the
schemes $\cX_i \subseteq \cT_{\fg_i}$. According to Proposition \ref{BPTS3}, we have
\[ (\ast) \ \ \ \ \ \ \ \ \ \ \ \ \ \  \varphi(\ft_i) \subseteq \fg_i \ \ \ \ \forall \ \varphi \in \cX(k) \ \ \text{and} \ \ g(\ft_i) = \ft_i \ \ \ \ \forall \ g \in S(\fg,\ft)\]
for each $i \in \{ 1,2 \}$. Consequently, Proposition \ref{BPTS3} shows that
\[ S(\fg,\ft) \lra S(\fg_1,\ft_1)\times S(\fg_2,\ft_2) \ \ ; \ \ g \mapsto (\pi_1\circ g \circ \iota_1, \pi_2\circ g \circ \iota_2)\]
is an injective homomorphism of groups.

Given $i \in \{1,2\}$, the map
\[ \zeta_i : \begin{cases} \cX_i(k) \lra \cX(k) \\ \varphi_i \mapsto \varphi_i\oplus \iota_{3-i}\end{cases}\]
is a morphism of affine varieties. As a result, given $g_i \in S(\fg_i,\ft_i)$, the elements $g_1\oplus \id_{\ft_2}$ and $\id_{\ft_1}\oplus g_2$ belong to $S(\fg,\ft)$, and the above map is
also surjective. \end{proof}

\bigskip

\section{Automorphisms and generic tori for Lie algebras of Cartan type} \label{S:AGT}
{\it Henceforth, we will assume that $k$ is an algebraically closed field of characteristic $p \geq 3$}. We shall apply the results of Section \ref{S:GTW} in the situation where $\fg$ is a restricted
Lie algebra of Cartan type. By definition, these restricted simple Lie algebras possess a {\it restricted $\ZZ$-grading}
\[ \fg = \bigoplus_{i=-r}^s \fg_i \ \ ; \ \ [\fg_i,\fg_j] \subseteq \fg_{i+j} \ \ ; \ \ \fg_i^{[p]} \subseteq \fg_{pi} \ \ \ \ \ (r,s\ge 1),\]
where we set $\fg_i := (0)$ for $i \not\in\{-r,\ldots,s\}$. Given such an algebra $\fg$, we consider the associated descending  filtration $(\fg_{(i)})_{i\ge -r}$, defined via
\[ \fg_{(i)} := \sum_{j\ge i} \fg_j.\]
By definition, the restricted Lie algebras of Cartan type are $p$-subalgebras of the algebra of derivations of the truncated polynomial rings
\[ \fA_n := k[X_1,\ldots, X_n]/(X_1^p,\ldots,X_n^p) \ \ \ \ \ \ (n\ge 1),\]
whose canonical generators will be denoted $x_1, \ldots, x_n$.

The Lie algebra $W(n) := {\rm Der}_k(\fA_n)$ is called the {\it $n$-th Jacobson-Witt algebra}. Its $p$-map is the standard $p$-th power of linear operators. We let $\partial_i \in W(n)$
denote the partial derivative with respect to the variable $x_i$. Then $\{\partial_1, \ldots,\partial_n\}$ is a basis of the $\fA_n$-module $W(n)$, so that $\dim_k W(n) = np^n$.

Here and further, we use the standard multi-index notation and write $x^\alpha :=x_1^{\alpha_1}\cdots x_n^{\alpha_n}$ and $|\alpha|:=\sum_{i=1}^n\alpha_i$ for
$\alpha=(\alpha_1,\ldots,\alpha_n)\in \NN_0^n$. Setting $\tau := (p\!-\!1,\ldots,p\!-\!1)$, the Lie algebra $W(n)$ obtains a restricted grading $W(n) = \bigoplus_{j=-1}^{|\tau|-1}W(n)_j$ via
\[ W(n)_j :=\sum_{i=1}^n\sum_{|\alpha|=j+1} kx^\alpha\partial_i.\]
We briefly describe the other Cartan type Lie algebras in terms of differential forms; the reader is also referred to \cite[Chap.~IV]{SF} and \cite{KS} for further details.

We set $\Omega^0(\fA_n)=\fA_n$ as well as $\Omega^1(\fA_n):=\Hom_{\fA_n}(W(n),\fA_n)$. The latter space carries canonical structures of an $\fA_n$-module and a $W(n)$-module
(cf.\ \cite{KS} or \cite{SF}) via
\[ (u.f)(D) := u.f(D) \ \ \text{and} \ \ (E.f)(D) = E(f(D))-f([E,D]),\]
for all $u \in \fA_n, \, D,E \in W(n), \, f \in \Omega^1(\fA_n)$. Let $\Omega^r(\fA_n):=\bigwedge^r \Omega^1(\fA_n)$ be the $r$-fold exterior power over $\fA_n$ and denote the exterior algebra by
\[\Omega(\fA_n):=\bigoplus_{r=0}^n \Omega^r(\fA_n).\]
The action of $W(n)$ on $\Omega^1(\fA_n)$ extends uniquely to an action by derivations on $\Omega(\fA_n)$. For $r \ge 2$, we have
\[ D.(f_1\wedge f_2\wedge \cdots \wedge f_r) := \sum_{i=1}^r f_1\wedge f_2\wedge \cdots \wedge f_{i-1} \wedge D.f_i\wedge f_{i+1}\wedge \cdots \wedge f_r,\]
for all $D \in W(n), \, f_i \in \Omega^1(\fA_n)$, so that $\Omega^r(\fA_n)$ is a $W(n)$-submodule.

Consider the $k$-linear map $\sfd: \Omega^0(\fA_n)\lra \Omega^1(\fA_n)$, given by
\[ \sfd u(D)=D(u) \ \ \ \ \ \forall \ u\in \fA_n, \ D\in W(n).\]
Associated with the {\it Cartan differential forms} $\omega_S=\sfd{x}_1\wedge\cdots \wedge \sfd{x}_n$; $\omega_H=\sum_{i=1}^r \sfd{x}_i\wedge \sfd{x}_{i+r}$ for $n=2r$, and
$\omega_K=\sfd{x}_{2r+1}+\sum_{i=1}^{r} (x_{i+r }\sfd{x}_i - x_{i }\sfd{x}_{i+r})$ for $n=2r\!+\!1$, there are the simple Cartan type Lie algebras $X(n)=X''(n)^{(2)}$, where
$X''(n)=\{D\in W(n) \ ; \  D.\omega_X=0\}$ if $X=S, H$ and $K''(n) := \{D \in  W(n) \ ; \ D.\omega_K \in \fA_n \omega_K\}$.

Recall that, given a $2r$-dimensional $k$-vector space with a non-degenerate symplectic form $\langle \, , \, \rangle : V\times V \lra k$, the group
\[ G\Sp(V) := \{ g \in \GL(V) \ ; \ \exists \ c(g) \in k^\times \ , \ \langle g(x),g(y)\rangle = c(g)\langle x,y\rangle \ \ \forall \ x,y \in V\}\]
is referred to as the {\it group of similitudes} of the symplectic space $(V,\langle\, ,\,\rangle)$.

\bigskip

\begin{Theorem}[\cite{Wi1,Wi2}]\label{AGT1} Let $\fg$ be a restricted Lie algebra of Cartan type $X(n)$ over a field $k$ of characteristic $p \geq 3$. If $\fg = W(1)$ or $H(2)$, assume that $p \geq 5$. Then the following statements hold:

{\rm (1)} \ The group $G := \Aut_p(\fg)$ is connected.

{\rm (2)} \ The group $G$ is a semidirect product $G= U\rtimes G_0$, where $G_0$ consists of those automorphisms preserving the $\ZZ$-grading of $\fg$, and $U:=\{g\in G \ ; \ (g\!-\!\id_\fg)(\fg_{(i)}) \subseteq \fg_{(i+1)}\}$. Furthermore, $U$ is unipotent and $G_0$ is reductive with
\[ G_0\cong \left \{ \begin{array}{ll} \GL(n), & \text{ for } X=W,S, \\
                                                              G \Sp(2r), & \text{ for } X=H \text{ and } K \text{ with } n=2r \text{ and } 2r\!+\!1, {\text{ respectively}}.
 \end{array} \right. \]

{\rm (3)} \ We have $g(\fg_{(i)}) = \fg_{(i)}$ for every $g \in G$ and $i \in \ZZ$.\end{Theorem}

\begin{proof} According to Wilson's results \cite[Thm. 2]{Wi1}, the theorem holds for the group of automorphisms which preserve the standard filtration. So we need to show that under the
assumptions of the theorem, all automorphisms of $\fg$ preserve it. Since this is the standard filtration defined with respect to the maximal subalgebra $\fg_{(0)} \subsetneq \fg$
\cite[Definition 3.5.1]{St}, we have to show that any automorphism preserves this subalgebra. When $p \ge 5$, this is Kreknin's theorem \cite{Kre}. For $p = 3$, the result is still true by \cite[(4.2.6)]{St},
except when $\fg = W(1)$, $H(2)$ or $K(3)$. When $\fg = K(3)$, \cite[(4.2.6)]{St} states that $\fg_{(0)}$ is the unique maximal subalgebra $\fh$ of
codimension 3 such that the $\fh$-module $\fg / \fh$ is not irreducible. It readily follows that $\fg_{(0)}$ is stable under all automorphisms of $\fg$ in all relevant cases, so the proof is complete. \end{proof}

\bigskip
\noindent
Owing to (3) of Theorem \ref{AGT1}, the subsets $\Tor(\fg)_{\fg_{(i)}}$ and $\Tor(\fg)^{\fg_{(i)}}$ of $\Tor(\fg)$ are $G$-invariant. It also follows that the function
\[ f_0 : \Tor(\fg)_{\fg_{(-1)}} \lra \NN_0 \ \ ; \ \ \ft \mapsto \dim_k \ft\cap\fg_{(0)}\]
is constant on the $G$-orbits of $\Tor(\fg)_{\fg_{(-1)}}$.

The following result, which states that $f_0$ forms a complete set of invariants for the $G$-orbits of $\Tor(\fg)_{\fg_{(-1)}}$, is based on Demu\v{s}kin's work \cite{De1,De2}, as corrected
in \cite[\S7]{St}. In our statement below, we represent the contact algebra as in \cite[(IV.5)]{SF}.

\bigskip

\begin{Theorem}\label{AGT2} Let $\fg$ be a Lie algebra of Cartan type $X(n)$, where $X\in \{W,S,H,K\}$. Then the following statements hold:

{\rm (1)} \ Two tori $\ft, \ft' \in \Tor(\fg)_{\fg_{(-1)}}$ belong to the same $G$-orbit if and only if $f_0(\ft) = f_0(\ft')$.

{\rm (2)} \ Each of the following tori represents the unique orbit $G\dact\ft_0$ of $\Tor(\fg)_{\fg_{(-1)}}$ with $f_0(G\dact\ft_0) = \{0\}$:

\indent \indent
{\rm (W)} \ $\ft_0=\langle (1\!+\!x_1)\partial_1,\ldots, (1\!+\!x_n)\partial_n\rangle$ for $\fg = W(n)$

\indent \indent
{\rm (S)} \ $\ft_0=\langle (1\!+\!x_1)\partial_1-x_n\partial_n,\ldots,(1\!+\!x_{n-1})\partial_{n-1}-x_n\partial_n\rangle$ for $\fg = S(n)$

\indent \indent
{\rm (H)} \ $\ft_0=\langle (1\!+\!x_1)\partial_1-x_{r+1}\partial_{r+1},\ldots,(1\!+\!x_r)\partial_{r+1}-x_{2r}\partial_{2r}\rangle$ for $\fg = H(2r)$

\indent \indent
{\rm (K)} \ $\ft_0=\langle x_1(1\!+\!x_{r+1}),\ldots, x_r(1\!+\!x_{2r}),\sum_{i=1}^rx_ix_{r+i}\!-\!x_{2r+1}\rangle$ for $\fg = K(2r\!+\!1)$.\hfill $\square$ \end{Theorem}

\bigskip
\noindent
The foregoing results yield the following description of the toral stabilizers of the Lie algebras of Cartan types $W,S$ and $H$.

\bigskip

\begin{Proposition} \label{AGT3} Let $\fg$ be a Lie algebra of Cartan type $W,S$ or $H$. Then $\ft_0$ is a generic torus and
\[ W(\fg,\ft_0) \cong S(\fg).\] \end{Proposition}

\begin{proof} According to Lemma \ref{GTWG7} and Theorem \ref{AGT1}, $\Tor(\fg)_{\fg_{(-1)}}$ is a closed, $G$-invariant subset of $\Tor(\fg)$. By the same token,
\[ \cO := \Tor(\fg)_{\fg_{(-1)}} \cap \Tor(\fg)^{\fg_{(0)}}\]
is an open, $G$-invariant subset of $\Tor(\fg)_{\fg_{(-1)}}$. In view of Theorem \ref{AGT2}(1), we thus have $\cO = G\dact\ft_0$. In particular, the orbit $G\dact\ft_0$ is open in
$\Tor(\fg)_{\fg_{(-1)}}$. Since $\fg = \fg_{(-1)}$ for $\fg = W(n), S(n), H(2r)$, it follows from the irreducibility of $\Tor(\fg)$ that $\ft_0$ actually is a generic torus of $\fg$. Our
second assertion is now a direct consequence of Theorem \ref{GTWG5}. \end{proof}

\bigskip

\begin{Remark} If $\fg$ is of type $K$, then $G.\ft_0$ is a subset of $\fg_{(-1)}$. Since $S_\fg \setminus \fg_{(-1)}$ is not empty, the $G$-saturation of $\ft_0$ is contained in the proper,
closed subset $\bar{S}_\fg\cap\fg_{(-1)} \subsetneq \bar{S}_\fg$, and Proposition \ref{GTWG4} shows that the torus $\ft_0$ is not generic. In fact, we will see in (\ref{NGT3}) that contact algebras do not afford generic tori. \end{Remark}

\bigskip

\section{Embeddings of $W(n)$} \label{S:EW}
In preparation for our computation of the toral stabilizers and Weyl groups of the generic tori of the Lie algebras of Cartan types $W, S$ and $H$, we provide in this section embeddings
$W(\mu(\fg)) \hookrightarrow \fg$ for such Lie algebras.

The map
\[ \Div : W(n) \lra \fA_n \ \ ; \ \ \sum_{j=1}^n f_j\partial_j \mapsto \sum_{j=1}^n\partial_j(f_j)\]
is called the {\it divergence} of $W(n)$. The special algebra $S(n)$ is the derived algebra of the subalgebra $\ker \Div$, cf.\ \cite[(IV.3)]{SF}.

Let $r \le n$. In the sequel, we shall identify $\fA_r$ with the subalgebra of $\fA_n$ that is generated by $x_1, \ldots, x_r$. A derivation $D \in {\rm Der}_k(\fA_r)$ will be regarded
as an element of ${\rm Der}_k(\fA_n)$ by setting $D(x_j) = 0$ for $r\!+\!1\le j \le n$.  As before, we shall employ multi-index notation. In particular, $\varepsilon_j$ denotes the vector, whose coordinates are given by the Kronecker symbols $\delta_{ij}$.

\bigskip

\begin{Lemma} \label{EW1} Let $n \ge 2$. Then
\[ \sigma_n : W(n\!-\!1) \lra S(n) \ \ ; \ \ D \mapsto D- \Div(D)x_n\partial_n\]
is an injective homomorphism of restricted Lie algebras.\end{Lemma}

\begin{proof} Let $D$ be an element of $W(n\!-\!1)$. Since
\[ \Div(\sigma_n(D)) = \Div(D) - \partial_n(\Div(D)x_n) = \Div(D) -\Div(D) = 0,\]
we have $\sigma_n(W(n\!-\!1)) \subseteq \ker \Div$. According to \cite[(IV.3.1)]{SF}, the map $\Div$ is a derivation from $W(n)$ into the module $\fA_n$. Given $D,E \in W(n\!-\!1)$, we thus obtain
\begin{eqnarray*}
[\sigma_n(D),\sigma_n(E)] & = & [D- \Div(D)x_n\partial_n,E- \Div(E)x_n\partial_n]\\
& = & [D,E] - D(\Div(E))x_n\partial_n+E(\Div(D))x_n\partial_n + [\Div(D)x_n\partial_n,\Div(E)x_n\partial_n] \\
& = & [D,E] - \Div([D,E])x_n\partial_n + [\Div(D)x_n\partial_n,\Div(E)x_n\partial_n]\\
& = &  \sigma_n([D,E]),
\end{eqnarray*}
so that $\sigma_n : W(n\!-\!1) \lra \ker \Div$ is a homomorphism of Lie algebras. As $W(n\!-\!1)$ is simple, we have
\[ \im \sigma_n = \sigma_n([W(n\!-\!1),W(n\!-\!1)]) \subseteq (\ker\Div)^{(1)} = S(n).\]
Moreover, the map $\sigma_n$ is obviously injective.

In view of Jacobson's formula, it suffices to verify the identity $\sigma_n(x^{[p]}) = \sigma_n(x)^{[p]}$ on a basis of $W(n\!-\!1)$. Let $a \le \tau$ be an element of $\NN_0^{n-1}$. By
virtue of \cite[(IV.2.7)]{SF}, we have
\[ (x^a\partial_j)^{[p]} = \left\{ \begin{array}{cl} 0 & a \ne \varepsilon_j \\ x^{\varepsilon_j}\partial_j & a = \varepsilon_j\end{array} \right.\]
for every $j \in \{1,\ldots,n\!-\!1\}$. Since
\[ \sigma_n(x^a\partial_j) = x^a\partial_j - a_jx^{a-\varepsilon_j}x_n\partial_n\]
there exist $\lambda_i \in k$ such that
\[\sigma_n(x^a\partial_j)^{[p]}(x_i) = \lambda_i \, x^{pa+\varepsilon_i-p\varepsilon_j} \ \ \ \ \ \ \forall \ i \in\{1,\ldots,n\!-\!1\}.\]
This expression vanishes unless $a= \varepsilon_j$. In the remaining case, we readily obtain $\sigma_n(x_j\partial_j) = x_j\partial_j - x_n\partial_n$, so that $\sigma_n(x_j\partial_j)^{[p]}
= \sigma_n(x_j\partial_j)$. \end{proof}

\bigskip
\noindent
We denote by $\cP(2r)$ the {\it Poisson algebra}, i.e., the space $\fA_{2r}$, endowed with the standard Poisson bracket $\{\,,\,\}$. By definition, $\{\,,\,\}$ is an associative bi-derivation,
defined on the generators $x_1,\ldots,x_{2r}$ by the rule:
\[\{ x_i , x_{i+r} \} = 1 = -\{x_{i+r},x_i\}  \ \ \ \ \ \ 1 \le i \le r,\]
and the remaining $\{x_i , x_j \} = 0$ (see \cite[p.168]{SF}). The center of $\cP(2r)$ consists of the space $k$ of constant polynomials, and there is an exact sequence
\[ (0) \to k \to \cP(2r) \stackrel{D_H}{\lra} H''(2r)\]
of Lie algebras, where $H''(2r)$ is the annihilator of the Hamiltonian form (see \cite[p.163]{SF}). In particular, $H''(2r)$ is a $p$-subalgebra of $W(2r)$. The image $H'(2r)$ of $D_H$ has
codimension 1 in $H''(2r)$, and the derived algebra $H'^{(1)} = H(2r)$ is the {\it Hamiltonian algebra}.

The Lie algebra $\cP(2r)$ is restrictable, and any $p$-mapping of $\cP(2r)$ renders the above exact sequence an exact sequence of restricted Lie algebras. We define a $p$-mapping by
the following formula (cf.\ \cite[p.~403]{St}):
\begin{equation} \label{pMappingPoisson}
(x^a)^{[p]} = \left\{ \begin{array}{ccl} 0 & \mbox{for} & a \not \in \{ 0 , \varepsilon_i + \varepsilon_{i+r} \ ; \ i = 1, \ldots , r \}, \\
x^a & \mbox{for} & a \in \{ 0 , \varepsilon_i + \varepsilon_{i+r} \, ; \,i = 1, \ldots , r \}. \end{array}\right.
\end{equation}
In particular, $1^{[p]} = 1$, and the center is a torus.

\bigskip

\begin{Lemma} \label{EWP} Let $r\ge 1$. Then the map
\[ \varphi_r : W(r) \lra \cP(2r) \ \ ; \ \ \sum_{j=1}^r f_j\partial_j \mapsto \sum_{j=1}^r x_j f_j(x_{r+1},\ldots,x_{2r})\]
is an injective homomorphism of restricted Lie algebras. \end{Lemma}

\begin{proof} We check that $\varphi_r$ is a Lie algebra homomorphism. Given $f,g \in \fA_r$ and $i,j \in \{1,\ldots,r\}$, we have $[f\partial_i,g\partial_j] = f\partial_i(g)\partial_j - g\partial_j(f)\partial_i$ in $W(r)$. Writing $f(x')$ for $f(x_{r+1},\ldots, x_{2r})$, we therefore obtain
\begin{eqnarray*}
\{ x_i f(x') , x_j g(x') \} & = & \{ x_i , x_j g(x') \} f(x') + x_i \{ f(x') , x_j g(x') \} \\
 & = & x_j \{ x_i , g(x') \} f(x') + x_i \{ f(x') , x_j  \} g(x') \\
 & = & x_j \partial_{i+r}(g(x')) f(x') - x_i \partial_{j+r} (f(x')) g(x')\\
 & = & \varphi_r([f\partial_i,g\partial_j]).
\end{eqnarray*}
Recall that $(x^a \partial_i )^{[p]} = 0$ if $a \neq \varepsilon_i$, and $(x_i \partial_i)^{[p]} = x_i \partial_i$ for all $i \in \{1, \ldots , r \}$. A comparison with
(\ref{pMappingPoisson}), shows that $\varphi_r \big( (x^a \partial_i )^{[p]} \big) =  \varphi_r(x^a \partial_i )^{[p]}$ for all $a \in \NN^r_0$. Consequently, $\varphi_r$ is an injective homomorphism of restricted Lie algebras. \end{proof}

\bigskip
\noindent
As a corollary, we readily obtain:

\bigskip

\begin{Lemma} \label{EW2} Let $r\ge 1$. The map $D_H \circ \varphi_r$ induces an embedding $W(r) \hookrightarrow H(2r)$ of restricted Lie algebras. \hfill $\square$ \end{Lemma}

\bigskip

\section{Toral stabilizers of Cartan type Lie algebras}\label{S:TSCA}
In this section, we turn to the computation of the toral stabilizers of restricted Lie algebras of Cartan type. Our approach involves embeddings of suitable $p$-subalgebras.

\subsection{Poisson algebras}

We begin by determining the toral stabilizers of the Poisson algebra with toral center. For $r \geq 1$, we consider the Lie algebra
\[ \fl_r := k z \oplus k x_1 \oplus \bigoplus_{i=1}^r k t_i,\]
whose bracket and $p$-mapping are given by
\[ [t_i , x_1] = \delta_{i1}x_1 \ \ \text{and} \ \ t_i^{[p]} = t_i \quad , \quad x_1^{[p]} = z = z^{[p]},\]
respectively, with all unspecified products being zero.

\bigskip

\begin{Lem}\label{Po4} For the restricted Lie algebra $(\fl_r , [p])$, we have $S(\fl_r) \simeq \ZZ / p\ZZ$. \end{Lem}

\begin{proof} Setting $\ft := \bigoplus_{i=2}^r kt_i$, we obtain a direct sum decomposition
\[ \fl_r = \fl_1\oplus \ft\]
 of $\fl_r$, whose second summand is a torus. Consequently, Corollary \ref{BPTS4} implies
\[S(\fl_r) \simeq S(\fl_1) \times S(\ft) \simeq S(\fl_1).\]
According to \cite[p.4215]{Fa2}, we have $S(\fl_1) \cong \ZZ/p\ZZ$, as desired. \end{proof}

\bigskip
\noindent
A restricted Lie algebra $\fg$ is referred to as {\it trigonalizable} if it can be embedded into a restricted Lie algebra of triangular matrices. This is equivalent to every restricted simple
$\fg$-module being one-dimensional.

\bigskip

\begin{Thm} \label{Po5} Let $\cP(2r)$ be the restricted Poisson algebra in $2r$ variables with toral center. Then $S(\cP(2r)) \simeq \GL_r(\FF_p)\ltimes (\FF_p)^r$. Moreover, for any
torus $\ft \subseteq \cP(2r)$ of maximal dimension, we have $S(\cP(2r),\ft) = \{ g \in \Aut_p(\ft)\ ; \ g(1) = 1 \}$. \end{Thm}

\begin{proof}  We first show that $S(\cP(2r)) \simeq \GL_r(\FF_p) \ltimes W$, where $W = (0)$ or $W=(\FF_p)^r$ with the natural action by $\GL_r(\FF_p)$. In view of Lemma \ref{EWP}, there exists an embedding $k \oplus W(r) \hookrightarrow \cP(2r)$ originating in the direct sum of the one-dimensional torus $k$ with $W(r)$. Let $\ft = k \oplus \ft'$, where $\ft' \subseteq W(r)$ is any torus of maximal dimension. Thanks to Lemma \ref{BPTS1}, Corollary \ref{BPTS4} and Proposition \ref{BPTS3}(3), we obtain inclusions
\begin{equation*}
\renewcommand{\arraystretch}{1.7}
\left( \begin{array}{c|c} 1  & 0 \\ \hline
 0 & S(W(r),\ft') \end{array} \right)
 \subseteq S(\cP(2r),\ft) \subseteq
\left( \begin{array}{c|c} 1  & (\FF_p)^r \\ \hline
 0 & \GL_r(\FF_p) \end{array} \right).
\end{equation*}
Owing to \cite[Thm.~1]{Pr2} and Proposition \ref{AGT3}, the group $S(W(r),\ft')$ is isomorphic to $\GL_r(\FF_p)$. Since the matrix group on the right identifies with the semi-direct product of $\GL_r(\FF_p)$ with the standard module $(\FF_p)^r$, we obtain embeddings
\[\GL_r(\FF_p) \ltimes (0) \subseteq S(\cP(2r),\ft) \subseteq \GL_r(\FF_p) \ltimes (\FF_p)^r.\]
Consequently, $S(\cP(2r),\ft) \simeq \GL_r(\FF_p) \ltimes W$, where $W := S(\cP(2r),\ft) \cap [(1)\! \ltimes\! (\FF_p)^r]$ is a $\GL_r(\FF_p)$-submodule of $(\FF_p)^r$. Since $(\FF_p)^r$ is irreducible, we have $W = (0)$ or $W=(\FF_p)^r$.

In order to show that $W \neq (0)$, we consider the subspace
\[\fl = k + k (1\! +\! x_{r+1}) +\sum_{i = 1}^r k x_i (1 \! + \! x_{i+r})\]
of $\cP(2r)$. We claim that $\fl$ is a $p$-subalgebra of $\cP(2r)$ which is isomorphic to $\fl_r$. Using the fact that the Poisson bracket is a bi-derivation for the associative structure, we obtain
\[\{ x_i(1\!+\!x_{i+r}), 1\!+\!x_{r+1}\} = \{x_i, 1\! +\! x_{r+1} \} (1\! +\! x_{i+r}) = \{ x_i , x_{r+1} \}(1\! +\! x_{i+r}) = \delta_{i1}(1\! +\! x_{i+r}).\]
Since polynomials in $x_i$ and $x_{i+r}$ commute with polynomials in the remaining variables, we see that all other brackets are $0$. Hence, $\fl \simeq \fl_r$ as abstract Lie algebras. We have $1^{[p]} = 1$ by definition. Also, using Jacobson's formula, we obtain
\[ (1 \! + \! x_{r+1})^{[p]} = 1^{[p]} + x_{r+1}^{[p]} = 1.\]
Setting $t_i = x_i (1\!+\!x_{i+r})$, we show that $t_i^{[p]} = t_i$. One can easily check $(\ad t_i)^p = \ad  t_i$ by comparing the action of these derivations on the associative
generators $x_1, \ldots, x_r$ and $1 \! + \! x_{r+1},\ldots,1 \! + \! x_{2r}$ of $\cP(2r)$. As the center $C(\cP(2r))$ coincides with $k$, there exists $\lambda \in k$ such that $t_i^{[p]} = t_i + \lambda 1$. Consider the decomposition $t_i = x_i \! + \! x_i x_{i+r}$. Jacobson's formula in conjunction with (\ref{pMappingPoisson}) yields
\[ t_i^{[p]} = x_i^{[p]} + (x_i x_{i+r})^{[p]} + \sum_{j = 1}^{p-1} s_j(x_i,x_i x_{i+r}) = x_i x_{i+r} + \sum_{j = 1}^{p-1} s_j(x_i , x_i x_{i+r}).\]
Observe that both, $x_i$ and $x_i x_{i+r}$, belong to the subspace $E := x_i \! \cdot \! k[x_{i+r}]$. It is easy to check that $E$ is stable under $\ad x_i$ and $\ad x_i x_{i+r}$. Hence all the $s_j(x_i,x_ix_{i+r})$ belong to $E$. Consequently, we obtain $\lambda 1 + t_i = t_i^{[p]} \equiv x_i x_{i+r} \mod E$, forcing $\lambda = 0$, as required.

According to \cite[(7.5.10)]{St}, the space $\ft := k + \sum_{i+1}^r k t_i \subseteq \fl$ is a torus of $\cP(2r)$ of maximal dimension. By Lemma \ref{Po4}, we have $S(\fl, \ft) \simeq \ZZ
/ p \ZZ$. Since the factor algebra $\fl / k$ is trigonalizable with torus $\ft': = \ft / k$ of maximal dimension, \cite[(6.1)]{Fa2} implies $S(\fl/k,\ft') = \{1\}$. Proposition \ref{BPTS3} now
provides an embedding $\xi$ in the top line of the following diagram:
\[ \begin{CD}
S(\fl,\ft) @>{\xi}>> \left( \begin{array}{c|c} 1 & \Lie_p(\ft',k) \\ \hline 0 & I_{\ft'} \end{array} \right) \\
@V{\alpha}VV @VVV \\
S(\cP(2r),\ft) @= \left( \begin{array}{c|c} 1 & W \\ \hline 0 & \GL_r(\FF_p) \end{array} \right).
\end{CD}\]
The map $\alpha$ is the inclusion given by Lemma \ref{BPTS1}, and the matrices on the right describe the endomorphisms of $\ft = k \oplus \ft'$ according to the block decomposition (\ref{IdentifyAutGroups}). By construction of all these maps, the diagram is commutative. Since $\im \xi \simeq S(\fl,\ft) \neq \{1\}$, we have $W \neq (0)$, as desired.

For the last statement of the theorem, observe that a maximal torus $\ft$ always contains $\fn = k.1$ as a subtorus. Corollary \ref{BPTS2} provides an inclusion $S(\cP(2r),\ft) \subseteq \{ g \in \Aut_p(\ft)\ ; \  g(1) = 1\}$. As these groups have same order, this inclusion is an equality. \end{proof}

\bigskip

\subsection{Lie algebras of Cartan type} The main result of this section reads:

\bigskip

\begin{Thm} \label{EW3} Let $\fg$ be a Lie algebra of Cartan type. Then there is an isomorphism
\[ S(\fg) \cong \GL_{\mu(\fg)}(\FF_p).\] \end{Thm}

\begin{proof} Assume first that $\fg$ is of type $W$. Proposition \ref{AGT3} provides a generic torus $\ft_0 \subseteq \fg$ and an
isomorphism
\[ W(\fg,\ft_0) \cong S(\fg).\]
In that case, the isomorphism $W(\fg,\ft_0) \cong \GL_{\mu(\fg)}(\FF_p)$ was established in \cite[Thm.~1]{Pr2}.

When $\fg$ is of type $S$ or $H$, Lemmas \ref{EW1} and \ref{EW2} furnish embeddings $W(\mu(\fg)) \hookrightarrow \fg$. Lemma \ref{BPTS1} thus yields $\GL_{\mu(\fg)}(\FF_p) \hookrightarrow S(\fg)$, so that $S(\fg) \cong \GL_{\mu(\fg)}(\FF_p)$.

The rest of the proof is devoted to contact algebras. Let $n=2r\!+\!1$ be an odd number. Following \cite[(IV.5)]{SF}, we endow the truncated polynomial ring $\fA_n$ with the structure of a
restricted Lie algebra; the corresponding Lie bracket $\langle \,  , \,  \rangle$ is called {\it contact bracket}. In the sequel, $K''(n)$ will denote the vector space $\fA_n$, endowed with this
bracket. Then $K(n) = K''(n)^{(1)}$ is the contact algebra.

Let $\ft := \bigoplus_{i = 1}^r k x_i x_{r+i} \oplus  k (1\!+\!x_n) \subseteq K(n)$. Then $\ft$ is a torus of maximal dimension in $K(n)$, and we claim that
\[ S(K(n),\ft) \simeq \GL_{r+1}(\FF_p). \]

We identify $\mathrm{Aut}_p(\ft)$ with $\GL_{r+1}(\FF_p)$, using the ordered basis $\{ 1\!+\!x_n, x_1 x_{r+1},\ldots,x_r x_{2r}\}$. Henceforth, all matrix representations will  describe automorphisms of $\ft$ by means of this basis.

By \cite[(7.5.15)]{St}, we know that the centralizer $\fc = C_{K''(n)}(1\!+\!x_n)$ is isomorphic to the Poisson algebra $\cP(2r)$ with toral center. Since $\ft \subseteq \fc$, Theorem \ref{Po5} yields the following matrix description:
\begin{equation} \label{lowerright}
\begin{pmatrix} 1 & (\FF_p)^r  \\  0 & \GL_r(\FF_p) \end{pmatrix} \equiv S(\fc,\ft) \subseteq S(K''(n),\ft),
\end{equation}
with the inclusion following from Lemma \ref{BPTS1}.

We will now embed lower triangular matrices into $S(K''(n),\ft)$, by means of a method similar to the one used in Theorem \ref{Po5}. Consider the subspace
\[  \fl := k x_r x_{2r} \oplus k x_r x_{2r} (1\!+\!x_n) \oplus k (1\!+\!x_n) \subseteq K''(n). \]
We first prove that $\fl$ is a restricted Lie algebra isomorphic to the one in \cite[p. 4216]{Fa2} (or to $\fl_1$ in Lemma \ref{Po4}). We check that such an isomorphism is defined by the assignment
\[ z \mapsto x_r x_{2r}, \quad x \mapsto x_r x_{2r} (1 \! + \! x_n), \quad t \mapsto 1 \! + \! x_n. \]
Indeed, directly from the rules given in \cite[p.173]{SF}, we obtain
\[ \langle x_r x_{2r},1 \! + \! x_n \rangle = 0 = \langle x_r x_{2r},x_r x_{2r}(1 \! + \! x_n)\rangle\]
as well as
\begin{eqnarray*}
\langle 1 \! + \! x_n,x_r x_{2r}(1 \! + \! x_n)\rangle & = & \langle 1, x_r x_{2r}(1 \! + \! x_n)\rangle + \langle x_n, x_r x_{2r}(1 \! + \! x_n)\rangle \\
                                                                & = & 2x_r x_{2r} + 2x_r x_{2r} x_n = 2x_r x_{2r}(1 \! + \! x_n).
\end{eqnarray*}
In view of \cite[p.177]{SF}, we also have
\[ (1 \! + \! x_n)^{[p]} = 1 + x_n^{[p]} = 1 \! + \! x_n.\]
Moreover, since $\langle x_r x_{2r},x_r x_{2r} x_n \rangle = 0$, we finally obtain
\[ (x_r x_{2r}(1 \! + \! x_n))^{[p]} = x_r x_{2r}.\]

Now let $\ft' := \sum_{i = 1}^{r-1}k x_i x_{i+r}$, so that $\ft = k (1 \! + \! x_n)+\ft' + k x_r x_{2r}$. Write $\ft_{\fl} = k(1 \! + \! x_n) + k x_r x_{2r}$ to ease notation. Since $\langle \ft' , \fl \rangle = 0$, there is an embedding $\ft' \oplus \fl \hookrightarrow K''(n)$ of restricted Lie algebras, which in turn yields an embedding
\[ S(\ft' \oplus \fl , \ft' \oplus \ft_{\fl} ) \simeq (1) \times S(\fl, \ft_{\fl}) \hookrightarrow S(K''(n),\ft)\]
of finite groups. According to \cite[p. 4216]{Fa2}, the group $S(\fl, \ft_{\fl})$ is generated by the transformation $x_r x_{2r} \mapsto x_r x_{2r}$, $(1 \! + \! x_n) \mapsto  x_r x_{2r} + (1 \! + \! x_n)$. Using the previous matrix identification, this yields
\[\begin{pmatrix} 1 & 0 \\ v & I_r \end{pmatrix}\in S(K''(n),\ft), \quad \text{where}\ v = \left[ \begin{matrix} 0 \\ \vdots \\ 0 \\ 1\end{matrix} \right].\]
In view of (\ref{lowerright}), the group $S(K''(n),\ft)$ also contains $\left(\begin{smallmatrix} 1 & 0 \\ 0 & \GL_r(\FF_p) \end{smallmatrix}\right)$, so the irreducibility of the $\GL_r(\FF_p)$-module $(\FF_p)^r$ yields
\begin{equation} \label{upperright} \begin{pmatrix} 1 & 0 \\ (\FF_p)^r & \GL_r(\FF_p) \end{pmatrix} \subseteq S(K''(n),\ft). \end{equation}
Using (\ref{lowerright}) and (\ref{upperright}), we see that $S(K''(n),\ft)$ contains the subgroup generated by all shear matrices, which, owing to \cite[Lemma 6.7.1]{Jac}, coincides with
$\SL_{r+1}(\FF_p)$. By (\ref{upperright}), it also contains matrices of arbitrary determinant, whence $S(K''(n),\ft) \cong \GL_{r+1}(\FF_p)$. Since the factor algebra $K''(n)/K(n)$ is
$p$-unipotent, Lemma \ref{BPTS1} now yields the result.
\end{proof}

\bigskip

\noindent Taking into account Theorem \ref{GTWG1} (1), we obtain the following remarkable corollary:

\bigskip

\begin{Cor} Let $\fg$ be a simple Lie algebra of Cartan type. Then the variety $\cT_\fg(k)$ is irreducible. \end{Cor}

\begin{proof}
Let $\ft \subseteq \fg$ be a torus of maximal dimension, so that $S(\fg,\ft) \subseteq \Aut_p(\ft)$. By Theorem \ref{GTWG1}, the variety $\cT_\fg(k)$ has $[\Aut_p(\ft)\!:\!S(\fg,\ft)]$ irreducible components. Now the previous theorem shows that $[\Aut_p(\ft)\!:\!S(\fg,\ft)] = 1$.
\end{proof}

\bigskip

\subsection{The Melikian algebra} For the sake of completeness, we mention the case of the restricted Melikian algebra. Recall that this is a simple graded Lie algebra over a field of
characteristic $p = 5$, which is neither classical nor of Cartan type. We refer the reader to \cite{Skr} for the details.

\bigskip

\begin{Thm} \label{EW4} Assume that $k$ has characteristic $p = 5$. Let $\cM$ be the restricted Melikian algebra. Then there is an isomorphism
\[ S(\cM) \cong \GL_{2}(\FF_p),\]
and the variety of embeddings $\cT_{\cM}(k)$ is irreducible. \end{Thm}

\begin{proof} By construction, there exists an embedding $W(2) \hookrightarrow \cM$ (see \cite[Sect. 1]{Skr}). By \cite[Cor. 4.4]{Skr}, all maximal tori in $\cM$ have dimension 2, so that
$\mu(W(2)) = \mu(\cM)$. Let $\ft \subseteq W(2) \subseteq \cM$ be a torus of dimension 2. Using Lemma \ref{BPTS1}, we obtain embeddings $S(W(2),\ft) \subseteq S(\cM,\ft)
\subseteq \Aut_p(\ft)$. By Theorem \ref{EW3}, this forces $S(\cM,\ft) = \Aut_p(\ft) \cong \GL_{2}(\FF_p)$. The fact that the variety of embeddings is irreducible now follows from
Theorem \ref{GTWG1}. \end{proof}

\bigskip

\section{Applications}\label{S:App}
The applications to be discussed in this final section pertain to weight space decompositions, the non-existence of generic tori and invariants for certain Lie algebras of Cartan type. While the second topic originates in Demu\v{s}kin's early work \cite{De1,De2} along with Strade's corrections \cite{St}, our observations on invariants are motivated by Premet's paper \cite{Pr2}.

\subsection{Irreducibility and Weights} Our first result illustrates why it is reasonable to consider toral stabilizers as a replacement for Weyl groups when dealing with arbitrary restricted Lie
algebras.

By definition, a {\it restricted $\fg$-module} $M$ is a $\fg$-module such that for every $x \in \fg$, the operator $m \mapsto x^{[p]}.m$ is the $p$-th power of the transformation of $M$
effected by $x$. If $\ft \subseteq \fg$ is a torus, then $M$ is a completely reducible $\ft$-module, giving rise to the {\it weight space decomposition}
\[ M = \bigoplus_{\lambda \in \Lambda_M} M_\lambda.\]
Each $\lambda \in \Lambda_M$ is a linear form satisfying $\lambda(t^{[p]})=\lambda(t)^p$ for every $t \in \ft$. As a result, $\Lambda_M \subseteq \Lie_p(\ft,k)$ is a subset of the
character group of $\ft$, i.e., the additive group of all homomorphisms $\ft \lra k$ of restricted Lie algebras, where the $p$-map on $k$ is the associative $p$-th power. The group
$\Aut_p(\ft)$ acts contragrediently on $\Lie_p(\ft,k)$ via
\[ h\dact \varphi := \varphi \circ h^{-1} \ \ \ \ \ \ \forall \ h \in \Aut_p(\ft), \ \varphi \in \Lie_p(\ft,k).\]
If $\ft \subseteq \fg$ is a torus of maximal dimension, then $S(\fg,\ft) \subseteq \Aut_p(\ft)$ also acts on $\Lie_p(\ft,k)$.

\bigskip

\begin{Thm} \label{IW1} Let $(\fg,[p])$ be a restricted Lie algebra, $\ft\subseteq \fg$ be a torus of maximal dimension. If $M$ is a restricted $\fg$-module with weight space decomposition
\[ M = \bigoplus_{\lambda \in \Lambda_M} M_\lambda,\]
then $\Lambda_M \subseteq \Lie_p(\ft,k)$ is $S(\fg,\ft)$-stable and the dimensions of the weight spaces are constant on the $S(\fg,\ft)$-orbits.\end{Thm}

\begin{proof} As in the proof of Lemma \ref{GTWG8}, put $A:=k[\cX_{\ft}]$ and let $j : \ft \hookrightarrow \fg\!\otimes_k\!A$ be the universal embedding of $\ft$. In particular, given a commutative $k$-algebra $R$, an element $\varphi_x: \ft \to \fg \otimes_k R$ corresponds to $x \in \Spec_k(A)(R)$ if and only if $\varphi_x = (\id_{\fg} \otimes x) \circ j$. Consider the restricted $A$-Lie algebra $\tilde{\fg}:=\fg\!\otimes_k A$. Then $\tilde{M} := M\!\otimes_k\!A$ is a restricted $\tilde{\fg}$-module, and the universal embedding gives rise to a weight space decomposition
\[ \tilde{M} = \bigoplus_{\gamma \in \Gamma_{\tilde{M}}}\tilde{M}_\gamma\]
of the $\ft$-module $\tilde{M}$, whose constituents are $A$-submodules of $\tilde{M}$. Let $x \in \Spec_k(A)(k)$ be a $k$-rational point. As observed in \cite[(4.2)]{Fa2}, we have
\[ M = \bigoplus_{\gamma \in \Gamma_{\tilde{M}}} M_\gamma(x),\]
where $M_\gamma(x) = M_{\gamma\circ\varphi_x^{-1}}$ is the weight space with weight $\gamma \circ \varphi_x^{-1}$ relative to the torus $\varphi_x(\ft) \subseteq \fg$. Since
each $\tilde{M}_\gamma$ is a finitely generated projective module over the integral domain $A$, it has constant rank and \cite[p.166]{Ma} yields
\[ \dim_kM_\gamma(x) = \rk(\tilde{M}_\gamma) \ \ \ \ \ \ \ \ \forall \ x \in \cX_\ft(k).\]
In particular, these dimensions are non-zero integers which depend only on $\gamma$.

Let $\iota : \ft \hookrightarrow \fg$ be the canonical inclusion. Given $g \in S(\fg,\ft)$, we have $g.\iota = \iota \circ g^{-1} \in \cX_{\ft}(k)$, and
\[ M = \bigoplus_{\gamma \in \Gamma_{\tilde{M}}} M_{\gamma\circ g}\]
is the weight space decomposition of $M$ relative to $\ft$. Thus, $\Lambda_M = g\dact\Gamma_{\tilde{M}}$ for all $g \in S(\fg,\ft)$ and $\dim_kM_{g\dact\lambda} = \rk(\tilde{M}_\lambda) = \dim_kM_\lambda$ for every $\lambda \in \Lambda_M$. \end{proof}

\bigskip

\begin{Cor} \label{IW2} Let $\ft \subseteq \fg$ be a torus of maximal dimension, $M = \bigoplus_{\lambda \in \Lambda_M}M_\lambda$ be a restricted $\fg$-module. If the variety
$\cT_{\fg}(k)$ is irreducible, then the following statements hold:

{\rm (1)} \ $\Lambda_M =\{0\}$, or $\Lambda_M$ contains $\Lie_p(\ft,k) \setminus \{0\}$.

{\rm (2)} \ All weight spaces of $M$ belonging to non-zero weights have the same dimension.\end{Cor}

\begin{proof} Let $\ft \subseteq \fg$ be a torus of maximal dimension, $M = \bigoplus_{\lambda \in \Lambda_M}M_{\lambda}$ be the corresponding weight space decomposition of $M$.
By Theorem \ref{IW1}, the set $\Lambda_M \subseteq \Lie_p(\ft,k)$ is $S(\fg,\ft)$-stable.

According to Theorem \ref{GTWG1}, the variety $\cT_{\fg}(k)$ has $[\Aut_p(\ft)\!:\!S(\fg,\ft)]$ irreducible components. In view of our assumption, this yields $S(\fg,\ft) = \Aut_p(\ft)$. Consequently $S(\fg,\ft) \cong \GL_{\mu(\fg)}(\FF_p)$ acts via two orbits on $\Lie_p(\ft,k) \cong \Hom_{\FF_p}(\FF_p^{\mu(\fg)},\FF_p)$, namely $\{0\}$ and $\Lie_p(\ft,k) \setminus\{0\}$. \end{proof}

\bigskip

\begin{Remark} Let $\fg$ be such that $\cT_{\fg}(k)$ is irreducible. For any restricted $\fg$-module $M = \bigoplus_{\lambda \in \Lambda_M}M_\lambda$, we have $\Lambda_M = \{0\}, \Lie_p(\ft,k)\setminus\{0\}$ or $\Lie_p(\ft,k)$. We consider the Jacobson-Witt algebra $\fg = W(n)$ acting via derivations in the truncated polynomial ring $\fA_n$. Then, $k \subseteq \fA_n$ is a submodule, and the three possibilities occur for $k$, $\fA_n / k$, and $\fA_n$, respectively. \end{Remark}

\bigskip
\noindent
We record the particular case given by the adjoint representation:

\bigskip

\begin{Cor} \label{IW3} Let $\fg = \fg_0 \oplus \bigoplus_{\alpha \in \Phi}\fg_\alpha$ be the root space decomposition of $\fg$ relative to a torus $\ft \subseteq \fg$ of maximal dimension. If the variety $\cT_{\fg}(k)$ is irreducible, then $\Phi \cup \{0\} = \Lie_p(\ft,k)$, and all root spaces of $\fg$ have the same dimension. \end{Cor}

\begin{proof} If $\Phi = \emptyset$, then $\fg = \fg_0$ is nilpotent and, being an irreducible variety of dimension zero, $\cT_\fg(k)$ is a singleton. Consequently, $\ft = (0)$, so that
$\Lie_p(\ft,k) = \{0\}$. The case where $\Phi \neq \emptyset$ follows directly from Corollary \ref{IW2}. \end{proof}

\subsection{Non-generic tori}
We will prove that the Weyl groups of the Poisson algebra $\cP(2r)$ do not coincide with the toral stabilizer, and obtain as a consequence that contact algebras do not afford generic tori.
This result strengthens Strade's observation \cite{St}, who showed that $K(2r\!+\!1)$ possesses infinitely many conjugacy classes of maximal tori.

To do so, we first gather some observations concerning automorphisms of $\cP(2r)$. Any Lie algebra automorphism of $\cP(2r)$ stabilizes the center $k \subseteq \cP(2r)$. The simple Lie
algebra $H(2r)$ can be realized as the derived subalgebra of the quotient algebra $\cP(2r) / k$. There results a homomorphism
\[ \pi : \Aut(\cP(2r)) \lra \Aut(H(2r))\]
of groups. Our analysis of $\Aut(\cP(2r))$ involves two classes of automorphisms. A Lie algebra automorphism $\sigma \in \Aut(\cP(2r))$ is a {\it semi-Poisson automorphism} if there exists
a non-zero scalar $\alpha \in k^\times$ such that $\alpha \sigma$ preserves the associative structure of $\cP(2r)$. If $\alpha = 1$, then $\sigma$ is referred to as a {\it Poisson
automorphism}. The set of semi-Poisson automorphisms is a subgroup of $\Aut(\cP(2r))$, which we denote by $\Aut_{\rm SP}(\cP(2r))$.

Let $\lambda \in \cP(2r)^*$ be a linear form which vanishes on $\cP(2r)^{(1)}$. Then
\begin{equation*}
\varphi_\lambda : \left\{ \begin{array}{ccl}
\cP(2r) & \to & \cP(2r) \\
f & \mapsto & f\!+\!\lambda(f)\end{array} \right.
\end{equation*}
is an automorphism of the Lie algebra $\cP(2r)$ and $\cH := \{ \varphi_\lambda \ ; \ \lambda \in (\cP(2r)/\cP(2r)^{(1)})^\ast\}$ is a subgroup of $\Aut(\cP(2r))$ such that $\Aut_{\rm
SP}(\cP(2r)) \cap \cH = \{\id_{\cP(2r)}\}$.

\bigskip

\begin{Thm}\label{NGT1} Let $\pi: \Aut(\cP(2r)) \to \Aut(H(2r))$ be the natural map. Then $\ker(\pi) = \cH$, and we have
\[ \Aut(\cP(2r)) = \Aut_{\rm SP}(\cP(2r)) \ltimes \cH.\]
\end{Thm}

\begin{proof} It follows from \cite[(7.3.6)]{St} that the restriction of $\pi$ to $\Aut_{SP}(\cP(2r))$ is surjective. We will proceed to verify $\ker(\pi) = \cH$. This will prove that the restriction of $\pi$ to $\Aut_{\rm SP}(\cP(2r))$ is also injective, so that $\Aut(\cP(2r)) = \Aut_{\rm SP}(\cP(2r)) \ltimes \cH$. Given $\sigma \in \ker(\pi)$, there exists a linear map
$\mu : \cP(2r)^{(1)} \to k$ such that $\sigma(f) = f + \mu(f)$ for all $f \in \cP(2r)^{(1)}$. Consequently,
\begin{eqnarray*}
\{ f , g \} & = & \{ f + \mu(f) , g + \mu(g) \} = \{ \sigma(f) , \sigma(g) \} \\
 & = & \sigma( \{ f , g \} ) = \{ f , g \} + \mu \big( \{ f , g \} \big),
\end{eqnarray*}
so that $\mu \big( \{ f , g \} \big) = 0$ for all $f,g \in \cP(2r)^{(1)}$. Since $\cP(2r)^{(1)}$ is a perfect Lie algebra, we have $\mu \equiv 0$, i.e.\ $\sigma(f) = f$ for every $f \in
\cP(2r)^{(1)}$.

We have a vector space decomposition $\cP(2r) = \cP(2r)^{(1)} \oplus k x^\tau$, where $x^\tau = x_1^{p-1} \cdots x_{2r}^{p-1}$ is the unique monomial of maximal degree $d =
2r(p\!-\!1)$. Let $\alpha \in k$ and $g \in \cP(2r)^{(1)}$ be such that $\sigma(x^\tau) = \alpha\, x^\tau + g$. Note that $g$ is a linear combination of monomials of degree lower than $d$. To see that $\sigma \in \cH$ we need to show that $\alpha = 1$ and $g \in k$.

Recall the notation $i' = i\!+\!r$ (resp. $i\! -\! r$) for $i \leq r$ (resp. $i > r$), so that $\{ x_i , x_j\} = \delta_{i' , j}$. Since the element
\[ \{ x_i , x^\tau \} = - x_1^{p-1} \cdots x_{i'}^{p-2} \cdots x_{2r}^{p-1},\]
lies in $\cP(2r)^{(1)}$, it is fixed by $\sigma$.  Hence,
\begin{eqnarray*}
- x_1^{p-1} \cdots x_{i'}^{p-2} \cdots x_{2r}^{p-1} & = & \sigma \big( \{ x_i , x^\tau \} \big) = \{ \sigma(x_i) , \sigma(x^\tau) \} = \{ x_i , \alpha\, x^\tau\!+\! g \} \\
& = & - \alpha\, x_1^{p-1} \cdots x_{i'}^{p-2}\cdots x_{2r}^{p-1} + \{ x_i , g \},
\end{eqnarray*}
so that
\[\{ x_i , g \} = (\alpha\! -\! 1) x_1^{p-1} \cdots x_{i'}^{p-2}\cdots x_{2r}^{p-1}\]
for all $i \in \{ 1 , \ldots , 2r \}$. Since bracketing by $x_i$ lowers the degree by $1$, $\{ x_i , g \}$ is a linear combination of monomials of degree $< d\!-\!1$. Consequently, both sides are
zero, so that $\alpha = 1$ and $\{ x_i , g \} = 0$ for all $i$. The latter condition implies that $g$ is central, i.e.\ $g \in k$ as desired. \end{proof}

\bigskip
\noindent
We turn to Weyl groups of the Poisson algebra.

\bigskip

\begin{Prop} \label{NGT2} Let $\ft \subseteq \cP(2r)$ be a torus of maximal dimension, $H \subseteq \Aut_p(\cP(2r))$ be a subgroup. Then the following statements hold:

{\rm (1)} \ There exists an embedding $\Nor_H(\ft) / \Cent_H(\ft) \hookrightarrow \FF_p^* \times \GL_r(\FF_p)$.

{\rm (2)} \ $W(\cP(2r),\ft) \subsetneq S(\cP(2r),\ft)$.

{\rm (3)} \ The Poisson algebra  $\cP(2r)$ affords no generic torus. \end{Prop}

\begin{proof} (1) Let $\fm \subseteq \cP(2r)$ denote the subspace of polynomials without constant term. Being the unique maximal ideal of the local algebra $\fA_{2r}$, $\fm$ is stable
under semi-Poisson automorphisms. Owing to \cite[(7.5.10)]{St}, the torus $\ft$ is conjugate under a Poisson automorphism to one of the following tori:
\[ T_q = k \oplus \sum_{i = 1}^q k (1\!+\!x_i)x_{i+r} \oplus \sum_{i=q+1}^r x_ix_{i+r}, \quad \ \ q \in \{ 0, \ldots,r \}.\]
(But the conjugating automorphism is {\bf not} restricted in general.) Each of these tori has the form $k \oplus T'$ for some torus $T' \subseteq \fm$. Since Poisson automorphisms
stabilize $\fm$, we also obtain a decomposition $\ft = k \oplus \ft_0$, where $\ft_0 \subseteq \fm$ is a torus.

Let $g \in \Nor_H(\ft)$. Thanks to Theorem \ref{NGT1}, the automorphism $g$ is a composite $g = g_0 \circ h$, where $g_0$ is a semi-Poisson automorphism and $h \in \cH$. Now, $\ft$
being a torus, we have $\ft \subseteq \cP(2r)^{(1)}$, so that $h$ acts trivially on $\ft$.  Thus, $g$ and $g_0$ induce the same automorphism of $\ft$, whence $g(k) \subseteq k$ and
$g(\ft_0) \subseteq \ft \cap \fm = \ft_0$. We obtain a well-defined mapping
\[ \Nor_H(\ft)\lra \Aut_p(k) \times \Aut_p(\ft_0)  \ \ ; \ \ g \mapsto ( g|_ k , g|_{\ft_0}),\]
which induces the desired embedding $\Nor_H(\ft) / \Cent_H(\ft) \hookrightarrow \FF_p^* \times \GL_r(\FF_p)$.

(2) According to Lemma \ref{GTWG2}, the Weyl group $W(\cP(2r),\ft)$ is contained in $S(\cP(2r),\ft)$. By part (1), there is an embedding $W(\cP(2r),\ft) \hookrightarrow \FF_p^{*} \times \GL_r(\FF_p)$, while Theorem \ref{Po5} yields $S(\cP(2r),\ft) \simeq \GL_r(\FF_p) \ltimes (\FF_p)^r$.

(3) This is a direct consequence of (2) and Theorem \ref{GTWG5}. \end{proof}

\bigskip

\begin{Thm} \label{NGT3} The contact algebra $K(2r\!+\!1)$ affords no generic torus. \end{Thm}

\begin{proof} We write $\fg := K''(2r\!+\!1)$ to ease notation, so that $\fg^{(1)} = K(2r\!+\!1)$. Note that $\fg^{(1)}$ is a $p$-ideal of $\fg$ such that the factor algebra is
$p$-unipotent. As a consequence, all tori of $\fg$ are contained in $\fg^{(1)}$, and we have $S_\fg = S_{\fg^{(1)}}$. Furthermore, if $\ft \subseteq \fg^{(1)}$ is a torus of maximal
dimension, then Lemma \ref{BPTS1} and Theorem \ref{EW3} yield $S(\fg,\ft) = S(\fg^{(1)},\ft) = \Aut_p(\ft)$.

Let us prove that restriction to the derived subalgebra induces a surjective homomorphism $\res: \Aut_p(\fg)^{\circ} \to \Aut_p(\fg^{(1)})$. By Theorem \ref{AGT1}, the group on the
right is connected. Thanks to \cite[Prop. 7.4 B]{Hu}, it therefore suffices to show the surjectivity of the map $\Aut_p(\fg) \to \Aut_p(\fg^{(1)})$. Setting $\Phi_\mu (D) := \mu \circ D\circ
\mu^{-1}$ for every $\mu \in \Aut_k(\fA_{2r+1})$ and $D \in W(2r\!+\!1)$, we obtain an isomorphism
\[ \Phi : \Aut_k(\fA_{2r+1}) \lra \Aut_p(W(2r\!+\!1)) \ \ ; \ \ \mu \mapsto \Phi_\mu\]
of groups. The automorphisms of $\fA_{2r+1}$ naturally act on differential forms, and we consider the subgroup
\[ G_K := \{ \mu \in \Aut_k(\fA_{2r+1}) \ ; \ \mu.\omega_K \in \fA_{2r+1}^{\times} \omega_K\}.\]
According to \cite[(7.3.2)]{St}, $\Phi$ induces an isomorphism $G_K \simeq \Aut_p(\fg^{(1)})$. From the definition $\fg = \{D \in  W(2r\!+\!1) \ ; \ D.\omega_K \in \fA_{2r+1}
\omega_K\}$, we see that for $\mu \in G_K$, the map $\Phi_{\mu}$ induces a restricted automorphism of $\fg$. Thus, any restricted automorphism of $\fg^{(1)}$ extends to $\fg$, as we
wanted to show.

In the sequel, we write $G := \Aut_p(\fg)^{\circ}$ and $G_K := \Aut_p(\fg^{(1)})$. Assume that there exists a generic torus $\ft \subseteq \fg^{(1)}$. In view of Proposition \ref{GTWG4}, the saturation $G_K.\ft$ is dense in $\bar{S}_{\fg^{(1)}}$. The above observations yield
\[ G.\ft = \res(G).\ft = G_K.\ft,\]
so that $G.\ft$ is a dense subset of $\bar{S}_\fg = \bar{S}_{\fg^{(1)}}$. Applying Proposition \ref{GTWG4} again, we conclude that $\ft$ is also a generic torus of $\fg$.

Note that the argument of the remark following Proposition \ref{AGT3} applies to any generic torus of $\fg^{(1)}$, yielding $\ft \not \subseteq \fg^{(1)}_{(-1)}$. Since $\ft \subseteq \fg^{(1)}$, this means $\ft \not \subseteq \fg_{(-1)}$. Owing to \cite[(7.5.14)]{St}, we may assume that $1\!+\!x_n \in \ft$. Then $\ft$ is also a torus of maximal dimension of the centralizer $\fc :=C_{\fg}(1\!+\!x_n)$. By \cite[(7.5.15)]{St}, there is an isomorphism $\varphi: \fc \to \cP(2r)$ of restricted Lie algebras such that $\varphi(1\!+\!x_n) = 1$, the canonical central element in the Poisson algebra. Applying Theorem \ref{Po5}, we obtain $S(\fc,\ft) = \Cent_{\Aut_p(\ft)}(1\!+\!x_n)$. Since $\ft \subseteq \fg$ is generic, Theorem \ref{GTWG5} in conjunction with the above observation $S(\fg,\ft) = \Aut_p(\ft)$ ensures that the map  $\Nor_{G}(\ft) \lra \Aut_p(\ft)$ is surjective. Thus, given $g \in \Aut_p(\ft)$, there exists $\varphi \in G$ such that $\varphi|_\ft = g$. Since $g$ fixes $1\!+\!x_n$, we have $\varphi|_{\fc} \in \Aut_p(\fc)$. Setting $H:=\Aut_p(\fc)$, we obtain a surjection $\Nor_{H}(\ft) \lra \Aut_p(\ft)$, which contradicts Proposition \ref{NGT2}(1). \end{proof}

\subsection{Invariants for Lie algebras of Cartan Type}\label{S:IL}
Given a restricted Lie algebra $(\fg,[p])$ with connected automorphism group $G = \Aut_p(\fg)^\circ$ and maximal torus $\ft \subseteq \fg$, we are interested in comparing the invariants
$k[\fg]^G$ and $k[\ft]^{S(\fg)}$ of the respective rings of polynomial functions. In case $S(\fg) = \GL_{\mu(\fg)}(\FF_p)$, the generators of the latter algebra were explicitly determined
by Dickson, cf.\ \cite{Di}.

We begin with some general observations concerning invariants of polynomial functions of a finite-dimensional associative $k$-algebra $A$. Let $G \subseteq \Aut(A)$ be a connected, closed
subgroup of the automorphism group of $A$. Given a finite-dimensional $A$-module $M$, we let $a_M : M \lra M$ be the left multiplication effected by the element $a \in A$. Let $T$ be an
indeterminate. For $a \in A$, we consider the characteristic polynomial
\[ P_M(T;a) := \det(T\!\cdot\!\id_M-a_M) \in k[T].\]
Then $\deg P_M(T;a) = \dim_kM =:n$ and the functions $\psi_{M,i} : A \lra k$, given by
\[ P_M(T;a) = \sum_{i=0}^n \psi_{M,i}(a)T^i\]
belong to the ring $k[A]$ of polynomial functions on $A$. Note that $\psi_{M,i} = \psi_{N,i}$ whenever $M \cong N$.

The group $G$ acts contragrediently on $k[A]$ via
\[ (g\dact\psi)(a) := \psi(g^{-1}(a)) \ \ \ \ \ \ \forall \ g \in G,\, \psi \in k[A],\, a \in A.\]
We record the following basic observation:

\bigskip

\begin{Lem} \label{IL1} Let $M$ be a finite-dimensional $A$-module. Then we have $\psi_{M,i} \in k[A]^G$ for every $i \in \{0,\ldots, \dim_kM\}$. \end{Lem}

\begin{proof} Given $g \in G$, we let $M^{(g)}$ be the $A$-module with underlying $k$-space $M$ and twisted action defined via
\[ a\dact m := g^{-1}(a).m \ \ \ \ \forall \ a \in A,\, m \in M.\]
Then we have $a_{M^{(g)}} = g^{-1}(a)_M$, so that
\[ \psi_{M^{(g)},i}(a) = \psi_{M,i}(g^{-1}(a))\]
for every $a \in A$ and $i \in \{0,\dots, \dim_kM\}$.

Let $S$ be a simple $A$-module. Since the connected group $G$ acts trivially on the set of primitive central idempotents of the algebra $ A/\Rad(A)$, it follows that
\[ S^{(g)} \cong S \ \ \ \ \forall \ g \in G.\]
Consequently,
\[ (g\dact \psi_{S,i})(a) = \psi_{S,i}(g^{-1}(a)) = \psi_{S^{(g)},i}(a) = \psi_{S,i}(a) \ \ \ \ \ \ \forall \ a \in A, \, g \in G,\]
whence $\psi_{S,i} \in k[A]^G$ for $i \in \{0,\ldots,\dim_kS\}$.

Let $(0) \lra M' \lra M \lra M'' \lra (0)$ be an exact sequence of finite-dimensional $A$-modules. The identity
\[ P_M(T;a) = P_{M'}(T;a)P_{M''}(T;a) \ \ \ \ \forall \ a \in A\]
implies that the $\psi_{M,i}$ belong to the subalgebra of $k[A]$ generated by the elements $\psi_{M',j}$ and $\psi_{M'',\ell}$. Our result now follows by induction on the length
of $M$. \end{proof}

\bigskip
\noindent
Returning to our standard set-up, we let $(\fg,[p])$ be a restricted Lie algebra with automorphism group $G = \Aut_p(\fg)^\circ$. Consider the two-sided ideal $I \unlhd U(\fg)$ of the
universal enveloping algebra that is generated by the set $\{x^p-x^{[p]} \ ; \ x \in \fg\}$. The finite-dimensional algebra
\[ U_0(\fg) := U(\fg)/I\]
is the {\it restricted enveloping algebra} of $\fg$, see \cite[(V.5.3)]{SF} for more details. By the universal property of $U_0(\fg)$, the group $G$ can be considered a connected subgroup of
the automorphism group of $U_0(\fg)$. Any finite-dimensional $U_0(\fg)$-module $M$ thus gives rise to invariant polynomials $\psi_{M,i}$, whose restrictions to $\fg$ are elements of the
ring $k[\fg]^G$ of invariants of polynomial functions on $\fg$. We shall henceforth consider the $\psi_{M,i}$ as elements of this ring.

Let $(\ft,[p])$ be a torus of dimension $m$. Then $\Aut_p(\ft) \cong \GL_m(\FF_p)$ acts contragrediently on the ring $k[\ft]$. The following subsidiary result interprets Dickson's work
\cite{Di} in this context.

\bigskip

\begin{Lem} \label{IL2} Let $(\ft,[p])$ be a torus of dimension $m$. Then
\[ k[\ft]^{\GL_m(\FF_p)} = k[\{\psi_{U_0(\ft),p^i} \ ; \ 0\le i \le m\!-\!1\}]\]
with algebraically independent generators. \end{Lem}

\begin{proof} Let $X(\ft)$ be the character group of $\ft$. By definition, $X(\ft)$ is the group of algebra homomophisms $\lambda : U_0(\ft) \lra k$. If $\{t_1,\ldots, t_m\}$ is a toral
basis of $\ft$, then the map $\lambda \mapsto (\lambda(t_1),\ldots, \lambda(t_m))$ provides an isomorphism $X(\ft) \cong \FF_p^m$ of groups. Since $\ft$ is a torus, the free
$U_0(\ft)$-module of rank $1$ decomposes into a direct sum of one-dimensional restricted $\ft$-modules:
\[ U_0(\ft) = \bigoplus_{\lambda \in X(\ft)} k_\lambda.\]
Consequently, we obtain for $t = \sum_{i=1}^m x_it_i \in \ft$ the identities
\[ P_{U_0(\ft)}(T;t) = \prod_{\lambda \in X(\ft)}(T\!-\!\lambda(t)) = \prod_{(a_1,\ldots,a_m)\in \FF_p^m} (T\!-\!a_1x_1\!-\cdots -\!a_mx_m).\]
Our result now follows directly from Dickson's Theorem \cite[(8.1.1)]{Di}. \end{proof}

\bigskip
\noindent
Let $M$ be a $U_0(\fg)$-module, $\ft \subseteq \fg$ be a torus of maximal dimension. Thanks to Theorem \ref{IW1}, the restrictions $\psi_{M,i}|_\ft$ belong to the subring $k[\ft]^{S(\fg)}$.

Given a commutative $k$-algebra $A$, we set $A^{p^r} := \{a^{p^r} \ ; \ a \in A\}$. The Krull dimension of $A$ will be denoted $\dim A$. Also, for a restricted Lie algebra $(\fg, [p])$, we define $d(\fg):= \max\{ \dim_k (k  x)_p\ ;\ x \in \fg \}$. By \cite[(8.6(3))]{Fa2}, we have $\mu(\fg) \leq d(\fg) \leq \rk(\fg)$. In particular, if $\fg$ admits a self-centralizing torus, then $\mu(\fg) = d(\fg) = \rk(\fg)$.

\bigskip

\begin{Lem} \label{IL3} Let $(\fg,[p])$ be a restricted Lie algebra and $\ft \subseteq \fg$ be a torus of maximal dimension. Set $\ell := d(\fg)\!-\!\mu(\fg)$. Then the following statements hold:

{\rm (1)} \cite[Theorem 2]{Pr3} There exists a $p$-polynomial $Q(T;x) \in k[\fg][T]$ of the form
\[ Q(T;x) = \sum_{i = 0}^{\mu(\fg)} \varphi_i(x) T^{p^{i+\ell}},\]
with $\varphi_{\mu(\fg)}=1$ and $\varphi_0 \neq 0$, such that for all $x \in \fg$, $\sum_{i = 0}^{\mu(\fg)} \varphi_i(x) x^{[p]^{i+\ell}} = 0$. Furthermore, all coefficients $\varphi_i$ belong to $k[\fg]^G$.

{\rm (2)} \ Let $M$ be a $U_0(\fg)$-module such that $M|_{U_0(\ft)}$ is free of rank $p^r$. Then $P_M(T;x) = Q(T;x)^{p^{r-\ell}}$.

{\rm (3)} \ For all $t \in \ft$, we have $Q(T;t) = P_{U_0(\ft)}(T;t)^{p^{\ell}}$.

{\rm (4)} \ For all $i \in \{ 0,\ldots,\mu(\fg)\}$, the function $x \mapsto \varphi_i(x)^{1/p^{\ell}}$ restricts to a rational function on $\bar{S}_\fg$.
\end{Lem}

\begin{proof} We set $\mu:=\mu(\fg)$ and $d:=d(\fg)$ to ease notation. Also, for $x \in \fg$ we denote by $x = x_s + x_n$ its Jordan-Chevalley decomposition, where $x_s$ and $x_n$ are $p$-semi-simple and $p$-nilpotent respectively, and $x_s, x_n \in (k x)_p$ (see  \cite[(II.3.5)]{SF}).

(1) Let $\nu$ be the smallest integer for which there exists a non-empty open subset $\Omega \subseteq \fg$ such that $\Omega^{[p]^\nu}\subseteq S_\fg$. By \cite[Theorem 2]{Pr3}, there exist polynomial functions $\varphi_0,\ldots,\varphi_{\mu} \in k[\fg]$, with $\varphi_\mu=1$, such that $\sum_{i=0}^{\mu} \varphi_i(x)\, x^{[p]^{\nu+i}} = 0$ for all $x \in \fg$. We need to check that $\nu = \ell$. Since $\varphi_\mu(x)=1$, we have $x^{[p]^{\nu+\mu}}\in \sum_{j=0}^{\nu+\mu-1}kx^{[p]^j}$, so that $\nu\!+\!\mu \le d$ and $\nu \le d\!-\!\mu=\ell$.

For the reverse inequality, we consider the subset $\Omega_0:=\{x \in \Omega\ ;\ \dim_k (kx)_p = d  \mbox{ and } \dim_k (kx_s)_p$ $ = \mu \}$ of $\fg$. We first show that it is not empty. It suffices to check that the following two sets are non-empty open subsets of $\fg$:
\[ \Omega':= \{x \in \fg\ ;\ \dim_k (kx)_p = d \} \quad \mbox{and} \quad \Omega'':= \{x \in \fg\ ;\ \dim_k (kx_s)_p = \mu \}. \]
By construction $\Omega'$ is not empty; for $\Omega''$ it follows for example from \cite[3.5(1)]{Fa2}. An element $x \in \fg$ belongs to $\Omega'$ if and only if the system of vectors $\{ x,x^{[p]},\ldots,x^{[p]^{d-1}} \}$ has maximal rank over $k$. This condition can be expressed in terms of non-vanishing of certain minors depending regularly on $x$, so $\Omega'$ is open. Regarding $\Omega''$, consider the morphism $q:\fg \to S_\fg \ ; x \mapsto x^{[p]^{\dim_k \fg}}$. The set $S_\fg^0:=\{x \in S_\fg \ ; \ \dim_k (k x)_p = \mu \}$ is an open subset of $S_\fg$, so that $\Omega''=q^{-1}(S_\fg^0)$ is open as well.

Given $x = x_s + x_n \in \Omega_0$, we have $x^{[p]^\nu} = x_s^{[p]^{\nu}}$. Our condition $\dim_k (kx_s)_p = \mu$ implies $x_s^{[p]^\mu} \in \sum_{i = 0}^{\mu - 1} x_s^{[p]^i}$, so that taking $p^\nu$-th powers yields
\[x^{[p]^{\nu + \mu}} = x_s^{[p]^{\nu + \mu}} \in \sum_{i = 0}^{\mu - 1} k\, x_s^{[p]^{\nu + i}} = \sum_{i = 0}^{\mu - 1} k\, x^{[p]^{\nu + i}}. \]
Consequently, $d = \dim_k (kx)_p \leq \nu\!+\!\mu$, whence $d\!-\!\mu \leq \nu$, as desired.

Finally, it follows from the construction in \cite{Pr3} that the polynomial $Q(T;x)$ is minimal for the annihilation property described in (1). The $G$-invariance of the coefficients $\varphi_i$ follows readily.

(2) First we show that $P_M(T;x)$ is a $p$-polynomial in $T$. Note that for any torus $\ft' \subseteq \fg$ of maximal dimension, the module $M|_{U_0(\ft')}$ is also free of rank $p^r$. Indeed, consider the weight space decompositions $M = \bigoplus_{\lambda \in \Lambda} M_{\lambda} =  \bigoplus_{\lambda' \in \Lambda'} M_{\lambda'}$ relative to $\ft$ and $\ft'$ respectively. By \cite[(4.2)]{Fa2}, there exists an isomorphism $\varphi:\ft' \to \ft$ such that $\lambda \mapsto \lambda \circ \varphi$ is a bijection $\Lambda \to \Lambda'$ satisfying $\dim_k M_{\lambda \circ \varphi} = \dim_k M_{\lambda}$ for every $\lambda \in \Lambda$. Because restricted modules over a torus are determined by weights and multiplicities of weight spaces, it readily follows that $M|_{U_0(\ft')}$ is free if and only if $M|_{U_0(\ft)}$ is free, of the same rank.

Let $\ft \subseteq \fg$ be a torus of dimension $\mu$. Since $M|_{U_0(\ft)} \cong U_0(\ft)^{p^r}$, we obtain
\[ P_M(T;t) = P_{U_0(\ft)}(T;t)^{p^r}\]
 for every $t \in \ft$. Let $\cO:=\{ t \in \ft\ ;\ \lambda(t) \neq \lambda'(t) \ \forall\, \lambda \neq \lambda' \in X(\ft) \}$, which is a dense open subset of $\ft$. For $t \in \cO$, the 
 polynomial $P_{U_0(\ft)}(T;t) \in k[T]$ is square-free, and its set of roots is an additive subgroup of $k$. Then \cite[\S 20]{Hu} implies that $P_{U_0(\ft)}(T;t)$ is a $p$-polynomial, 
 whence $\psi_{U_0(\ft_0),i} = 0$ for $i \not \in \{p^j \ ; \ 0 \le j \le \mu\}$.  Given $t \in \ft$, we obtain
\[ \psi_{M,i}(t) = \left\{ \begin{array}{cl} \psi_{U_0(\ft_0),p^j}^{p^r}(t) & \text{for}\ i = p^{r+j}, j \in \{0,\ldots, \mu\} \\ 0 & \text{otherwise.} \end{array} \right.\]

Suppose that $i$ is not of the form $p^{r+j}$ with $j \in \{0,\ldots, \mu\}$. Then the polynomial function $\psi_{M,i}(x)$ vanishes on $\fT := \bigcup_{\ft \in \Tor(\fg)} \ft$. By Lemma \ref{GTWG8}, $\fT$ is dense in $\bar{S}_\fg$, so that $\psi_{M,i}$ vanishes on $\bar{S}_\fg$. Thus we obtain
\begin{equation} \label{ppolynomiality} 
P_M(T;x) = \sum_{j = 0}^{\mu} \psi_{M,p^{r+j}}(x) T^{p^{r+j}} \ \ \ \ \forall \  x \in \bar{S}_\fg.
\end{equation}

For arbitrary $x=x_s+x_n \in \fg$, the elements $x_s$ and $x_n$ act as commuting semi-simple and nilpotent endomorphisms on $M$, so that $P_M(T;x) = P_M(T;x_s)$. Since $x_s \in \bar{S}_\fg$, the above formula implies $\psi_{M,i}(x) = 0$ unless $i = p^{r+j}$ for some $j \in \{ 0,\ldots,\mu\}$, showing that Formula (\ref{ppolynomiality}) holds for all $x \in \fg$.

Next we check that $P_M(T;x) = Q(T,x)^{p^{r-\ell}}$, considered as elements of $k[T]$. It suffices to verify the identity for all $x$ in some dense open subset of $\fg$. Since $P_M(T;x)$ and $Q(T;x)^{p^{r-\ell}}$ are monic $p$-polynomials in $T$ of the same degree, it is enough to show that they have the same set of roots. For each $x \in \fg$, denote by $\Rac_Q(x)$ and $\Rac_P(x)$ the set of roots of $Q(T;x)$ and $P_M(T;x)$. Since $Q(T;x)$ annihilates $x$ in $\fg$, it also annihilates $x_M \in \End_k(M)$, the action of $x$ in $M$. Therefore, it must be divisible by the minimal polynomial of $x_M$. Since the latter has the same roots as $P_M(T;x)$, we obtain $\Rac_P(x) \subseteq \Rac_Q(x)$ for all $x \in \fg$. Conversely, for $x = x_s + x_n$, we have $\Rac_P(x) = \Rac_P(x_s)$. When $x \in \Omega_0$ as above, the torus $(k  x_s)_p$ has dimension $\mu$, and $\Rac_P(x_s)$ is an elementary abelian $p$-group of rank $\mu$. On the other hand, since $Q(T;x)$ has degree $p^{d}$ and is the $p^{\ell}$-th power of a $p$-polynomial of degree $p^\mu$, $\Rac_Q(x)$ is an elementary abelian $p$-group of rank at most $d\!-\!\ell = \mu$. This forces  $\Rac_P(x) = \Rac_Q(x)$ when $x \in \Omega_0$, as we needed.

(3) readily follows from the proof of (2).

(4) Let $\beta$ be an alternating $\mu$-linear form on $\fg$, and let $\Omega_\beta:=\{x \in \fg\ ; \ \beta(x^{[p]^{\ell}},\ldots,x^{[p]^{d-1}}) \neq 0\}$. The collection of all $\Omega_\beta$ forms an open covering of the set of all $x \in \fg$ such that $\{ x^{[p]^{\ell}},\ldots,x^{[p]^{d-1}}\}$ is a linearly independent family. Using the relation $x^{[p]^{d}}=-\sum_{i=0}^\mu \varphi_i(x) x^{[p]^{\ell + i}}$, we obtain for all $i \in \{0,\ldots,\mu\!-\!1\}$:
\begin{equation} \label{expressionwithbeta}
 \varphi_i(x) = -\frac{ \beta( x^{[p]^{\ell}} , \ldots , x^{[p]^{d}} , \ldots , x^{[p]^{d-1}} )}{\beta( x^{[p]^{\ell}} , \ldots , x^{[p]^{\ell + i}} , \ldots , x^{[p]^{d-1}} ) } \ \ \ \ \forall \ x \in \Omega_\beta.
\end{equation}
Now let $\ft \in \Tor(\fg)$. If $\ft \cap \Omega_\beta \neq \emptyset$, then $\beta$ does not vanish on $\ft$, and there exists $\lambda \in k^\times$ such that $\beta(t_1^{[p]^{\ell}} , \ldots , t_\mu^{[p]^{\ell}}) = \lambda \beta(t_1,\ldots,t_\mu)^{p^{\ell}}$ for all $t_1,\ldots,t_\mu \in \ft$. Indeed, the $p^{\ell}$-th root of the left hand side defines an alternating $\mu$-linear form on $\ft$, and so must be proportional to $\beta$. It follows that
\begin{equation}
\varphi_i(x) = - \frac{\beta(x,\ldots,x^{[p]^{\mu}},\ldots,x^{[p]^{\mu-1}})^{p^\ell}}{\beta(x,\ldots, x^{[p]^{i}},\ldots,x^{[p]^{\mu-1}})^{p^\ell}} \ \ \ \ \ \forall \ x \in \ft \cap \Omega_\beta.
\end{equation}
This shows that $x \mapsto \varphi_i(x)^{1/{p^\ell}}$ is a rational function on the Zariski closure of $\fT \cap \Omega_\beta$. For a suitable choice of $\beta$, this is a non-empty open subset of $\fT$, so by Lemma \ref{GTWG8} its Zariski closure is $\bar{S}_\fg$.
\end{proof}

\bigskip
\noindent
We record the following corollary:
\bigskip

\begin{Cor} \label{IL6}
Let $(\fg,[p])$ be a restricted Lie algebra and $\ell = d(\fg)\!-\!\mu(\fg)$. Then the morphism $\fg \to \fg\ ;\ x \mapsto x^{[p]^\ell}$ induces a dominant morphism $\fg \to \bar{S}_\fg$. In particular, $\ell = 0$ if and only if $\mu(\fg) = \rk(\fg)$.
\end{Cor}

\begin{proof}
The first part is a rephrasing of the equality $\ell = \nu$ (in the notation of \cite{Pr3}), proved in Part (1) above. If $\ell = 0$ then $\fg = \bar{S}_\fg$, and \cite[(3.7)]{Fa2} yields $\rk(\fg)\!-\!\mu(\fg) = 0$. The reverse implication is clear.
\end{proof}

\bigskip
\noindent
Let $f : A \lra B$ be an injective homomorphism of finitely generated, commutative $k$-algebras. If there exists $r \in \NN_0$ such that $B^{p^r}\subseteq \im f$, then $f$ is called
an {\it inseparable isogeny of exponent $\le r$}. In that case, the comorphism $f^\ast : {\rm Maxspec}(B) \lra {\rm Maxspec}(A)$ is a homeomorphism.

If $G$ is a group acting on $A$ via automorphisms such that $A^G$ is finitely generated, then the variety ${\rm Maxpec}(A)/\!\!/ G := {\rm Maxspec}(A^G)$ is called the {\it algebraic quotient} of ${\rm Maxspec}(A)$ by $G$.  If, in addition, all orbits are closed, then ${\rm Maxpec}(A)/G:={\rm Maxspec}(A^G)$ is referred to as the {\it geometric quotient} of ${\rm
Maxspec}(A)$ by $G$, see \cite[(II.3.2)]{Kr} for more details.

\bigskip

\begin{Thm} \label{IL4} Let $(\fg,[p])$ be a restricted Lie algebra with connected automorphism group $G = \Aut_p(\fg)^\circ$ and $\ft$ a torus of maximal dimension. Set $\ell := d(\fg)\!-\!\mu(\fg)$.

{\rm (1)} \ The canonical restriction map $k[\fg] \lra k[\ft]$ induces an algebra homomorphism $\res  : k[\bar{S}_\fg]^G \to k[\ft]^{W(\fg,\ft)}$, whose image contains
$(k[\ft]^{\GL_{\mu(\fg)}(\FF_p)})^{p^\ell}$.

{\rm (2)} \ We have $\dim k[\bar{S}_\fg]^G \geq \mu(\fg)$. Equality holds if and only if ${\rm res}$ is injective.

\smallskip
\noindent
Assume further that $\fg$ affords a generic torus $\ft_0$. Then:

{\rm (3)}\ The algebra $k[\bar{S}_\fg]^G$ is finitely generated and of Krull dimension $\mu(\fg)$.

{\rm (4)} \ If $W(\fg,\ft_0) = \GL_{\mu(\fg)}(\FF_p)$, then ${\rm res} : k[\bar{S}_\fg]^G \hookrightarrow k[\ft_0]^{W(\fg,\ft_0)}$ is an inseparable isogeny of exponent $\le \ell$. In particular, the varieties $\bar{S}_\fg/\!\!/G$ and $\ft_0/W(\fg,\ft_0)$ are homeomorphic.
\end{Thm}

\begin{proof} We will write $\GL_\mu$ instead of $\GL_{\mu(\fg)}(\FF_p)$ to ease notation.

(1) Clearly, the restriction map induces a homomorphism $\res  : k[\bar{S}_\fg]^G \to k[\ft]^{W(\fg,\ft)}$. Let $Q(T;x)$ be as in Lemma \ref{IL3}, so that its coefficients $\varphi_i$ belong to $k[\fg]^G$. By Lemma \ref{IL3}(3), we have $\res (\varphi_i) = \psi_{U_0(\ft),p^i}^{p^{\ell}}$ for all $i \in \{0,\ldots,\mu\}$. In view of Lemma \ref{IL2}, the $\psi_{U_0(\ft),p^i}$ generate $k[\ft]^{\GL_\mu}$, whence $(k[\ft]^{\GL_\mu})^{p^\ell} \subseteq \im \res $.

(2) We have $(k[\ft]^{\GL_\mu})^{p^\ell} \subseteq A:=\im {\rm \res} \subseteq k[\ft]$. By the Theorem of Hilbert-Noether \cite[(1.3.1)]{Be}, $k[\ft]$ is a finitely generated module over the noetherian ring $k[\ft]^{\GL_\mu}$. So it is also a finitely generated module over $B:=(k[\ft]^{\GL_\mu})^{p^\ell}$. Being a $B$-submodule of $k[\ft]$, the algebra $A$ is noetherian as well. Consequently, $A$ is a finitely generated integral domain.

Owing to Lemma \ref{IL2}, the algebras $(k[\ft]^{\GL_\mu})^{p^\ell}$ and $k[\ft]$ have common Krull dimension $\mu(\fg)$. Using \cite[(8.2.1.A)]{Ei}, we obtain $\dim A = \mu(\fg)$. It readily follows that $\dim k[\bar{S}_\fg]^G \geq \mu(\fg)$. Furthermore, since $\bar{S}_\fg$ is irreducible (see \cite[\S 2]{Pr1} or \cite[(3.7)]{Fa2}), $k[\bar{S}_\fg]^G$ is an integral domain. Therefore, $\dim k[\bar{S}_\fg]^G = \dim A$ if and only if $\ker {\res} = (0)$.

(3) Owing to Proposition \ref{GTWG4}, the $G$-saturation $G.\ft_0$ lies dense in $\bar{S}_\fg$. This readily implies that $\res  : k[\bar{S}_\fg]^G \to k[\ft]^{W(\fg,\ft_0)}$ is injective, hence $k[\bar{S}_\fg]^G \simeq A$ is a finitely generated algebra of Krull dimension $\mu(\fg)$.

(4) readily follows from properties (1) to (3).
\end{proof}

\bigskip
\noindent
The following result extends Premet's Theorem \cite[Thm.~1]{Pr2} concerning $W(n)$, where $\bar{S}_{W(n)}
= W(n)$.

\bigskip

\begin{Thm} \label{IL5} Suppose that $p\ge 3$ and let $(\fg,[p])$ be a restricted Lie algebra of Cartan type $W,S$ or $H$ with generic torus $\ft \subseteq \fg$. Then the restriction map induces an isomorphism $k[\bar{S}_\fg]^G \stackrel{\sim}{\lra} k[\ft]^{\GL_{\mu(\fg)}(\FF_p)}$. In particular, $k[\bar{S}_\fg]^G$ is a polynomial ring in $\mu(\fg)$ variables.
\end{Thm}

\begin{proof} First we note that for $\fg = W(n), S(n)$ or $H(n)$, there exists a generic torus $\ft \subseteq \fg$ (Proposition \ref{AGT3}), and the corresponding Weyl group $W(\fg,\ft) \simeq \GL_{\mu(\fg)})(\FF_p)$ (Theorems \ref{GTWG5} and \ref{EW3}).

(1) Let $\fg=W(n)$. For $p=3$ and $\fg=W(1)$ our assertion follows from the Chevalley Restriction Theorem. Alternatively, consider the generic torus $\ft_0 \subseteq \fg$, so that $W(\fg,\ft) = \GL_n(\FF_p)$. By virtue of \cite[(IV.2.5)]{SF}, the algebra $W(n)$ has a self-centralizing maximal torus; in particular $\ell = d(\fg) - \mu(\fg)=0$. Now Corollary \ref{IL6} yields $\bar{S}_\fg = \fg$, and our assertion follows from Theorem \ref{IL4}(4).

(2) Let $\fg=H(2r)$ and consider the $W(2r)$-module $\fA_{2r}$. Let $\ft_0$ be a generic torus of $W(2r)$ such that $W(2r)=\ft_0\oplus W(2r)_{(0)}$. Since $\fA_{2r} \cong U_0(W(2r))\!\otimes_{U_0(W(2r)_{(0)})}\!k$ is a free $U_0(\ft_0)$-module of rank $1$, Lemma \ref{IL3} implies
\begin{equation} \label{gencharpolW}
P_{\fA_{2r}}(T;x) = \sum_{i = 0}^{2r} \psi_i(x)T^{p^i} \ \ \ \  \forall \ x \in W(2r),
\end{equation}
with each $\psi_i$ being a homogeneous polynomial function of degree $p^{2r}\!-\!p^i$ on $W(2r)$. Also note that, for any subtorus $\ft \subseteq \ft_0$ of dimension $r$, $\fA_{2r}$ is a free $U_0(\ft)$-module of rank $p^r$.

Let $G = \Aut_p(\fg)^\circ$. Then $P_{\fA_{2r}}(T;x)$ is the characteristic polynomial of $x \in \fg$ acting on $\fA_{2r}$. By Lemma \ref{IL1}, the functions $\psi_i|_\fg$ are $G$-invariant. We show that they are $p^r$-th powers in $k[\fg]^G$. Consider the Poisson algebra $\cP(2r)$ with $p$-unipotent center. For any $f \in \cP(2r)$, the linear map $\ccD_f = \{ f , - \}$ is a derivation of the associative algebra $\fA_{2r}$, so that $\ccD_f \in W(2r)$. According to \cite[Lemma 6.3]{Skr2}, there exist polynomial functions $\phi_0,\ldots,\phi_{r} \in k[\cP(2r)]$ such that
\begin{equation} \label{gencharpolP}
P_{\fA_{2r}}(T;\ccD_f) = \sum_{i = 0}^{r} \phi_i(f)^{p^r} T^{p^{r+i}} \ \ \ \ \forall \ f \in \cP(2r).
\end{equation}
The map $\ccD : \cP(2r) \lra W(2r) \ ;  \ f \mapsto \ccD_f$ induces a surjection $\cP(2r)^{(1)} \twoheadrightarrow \fg$. Hence there exists a linear map $\delta : \fg \lra \cP(2r)^{(1)}$
such that $\ccD\circ \delta = \id_\fg$. In view of (\ref{gencharpolW}) and (\ref{gencharpolP}), the polynomial maps $\varphi_i := \phi_i\circ\delta \in k[\fg]$ satisfy the identities
\[ \varphi_i(h)^{p^r} = \phi_i(\delta(h))^{p^r} = \psi_i(\ccD_{\delta(h)}) = \psi_i(h) \ \ \ \ \forall \ h \in \fg, \, i \in\{0,\ldots,r\}.\]
Since $\varphi_i^{p^r} = \psi_i$ is $G$-invariant, we readily obtain that $\varphi_i \in k[\fg]^G$ as well.

Let $\ft \subseteq \fg$ be a generic torus. Theorem \ref{IL4} shows that the restriction map induces an embedding $\res  : k[\bar{S}_\fg]^G \hookrightarrow k[\ft]^{\GL_r(\FF_p)}$. As noted before, $\fA_{2r}$ is a free $U_0(\ft)$-module of rank $p^r$, so
\[P_{\fA_{2r}}(T;t) = P_{U_0(\ft_0)}(T;t)^{p^r} = (\sum_{i=0}^r \kappa_i(t)T^{p^i})^{p^r} \quad \ \forall \ t \in \ft,\]
where the coefficients $\kappa_i \in k[\ft]$ are the Dickson generators of $k[\ft]^{\GL_{r}(\FF_p)}$. In view of the above, we obtain $\varphi_i(t) = \kappa_i(t)$ for all $i \in \{0,\ldots,r\!-\!1\}$ and $t \in \ft$, so that $\kappa_i = {\rm \res}(\varphi_i)$. Therefore, ${\rm \res}$ is surjective.

(3) Now we consider $\fg = S(n)$ and check that $\ell = d(\fg)\!-\!\mu(\fg) = 1$. If $\ell = 0$, then Corollary \ref{IL6} implies $\rk(\fg) = \mu(\fg)$, a contradiction. Since $S(n) \subseteq W(n)$, we have $d(\fg) \leq d(W(n)) = n$, so that $\mu(S(n))=n\!-\!1$ implies $\ell = d(\fg)\!-\!\mu(\fg) \leq 1$. By Theorem \ref{IL4} we have an embedding $\res :k[\bar{S}_\fg]^G \hookrightarrow k[\ft]^{\GL_{n-1}(\FF_p)}$. Consider the $p$-polynomial $Q(T;x) = \sum_{i = 0}^{n-1} \varphi_i(x) T^{p^{i+1}}\in k[\fg][T]$ given by Lemma \ref{IL3}; we know that the functions $\fg \to k\ ; \ x \mapsto \varphi_i(x)^{1/p}$ restrict to the Dickson generators of $k[\ft]^{\GL_{n-1}(\FF_p)}$. We will show that these functions are polynomial on $\fg$, which will prove that $\res$ is surjective.

We use the notation of Section \ref{S:AGT} for truncated polynomials and their derivations. Recall the decomposition $\fg = \fg_{(0)} \oplus \fg_{-1}$, where $\fg_{(0)}$ is a $G$-invariant subalgebra of codimension $n$. For every $x \in \fg$, there is a unique decomposition $(kx)_p = \ft_x \oplus \fn_x$, where $\ft_x$ is a torus and $\fn_x$ a $p$-unipotent subalgebra. We 
consider the set $\Omega:=\{ x \in \fg\ ;\ \dim_k (kx)_p = n,\ \dim_k \ft_x = n\!-\! 1,\ \ft \cap \fg_{(0)} = (0) \}$. We check that $\Omega$ is a non-empty open set: note that $\Omega=\Omega' \cap \Omega''$, where 
\[ \Omega':=\{ x \in \fg\ ;\ \dim_k (kx)_p = n \} \quad \mbox{and} \quad \Omega'':=\{ x \in \fg\ ;\  \dim_k \ft_x = n\!-\! 1,\ \ft_x \cap \fg_{(0)} = (0) \}. \]
Since $d(\fg)=\ell\!+\!\mu(\fg)=n$, the set $\Omega'$ is open and non-empty. For $\Omega''$, we have seen in the proof of Lemma \ref{IL3}(1) that the condition $\dim_k \ft_x = 
n\!-\!1$ defines a non-empty open subset of $\fg$. Given any such element $x$, the subspace $\ft_x$ is spanned by the linearly independent family $\{ x^{[p]^{\ell+1}}, \ldots, x^{[p]^{\ell+n-1}}\}$. The condition that $\ft_x \cap \fg_{(0)} = (0)$ is equivalent to the fact that the residue classes of $\{ x^{[p]^{\ell}}, \ldots, x^{[p]^{\ell+n-1}}\} $ in $\fg / \fg_{(0)}$ form a system of maximal rank, which is an open condition. Hence $\Omega''$ is open. We check that it is not empty. By Theorem \ref{AGT2} there exists a torus  $\ft_0$ of maximal dimension such that $\ft_0 \cap \fg_{(0)} =(0)$. By \cite[3.5(1)]{Fa2}, there exists $t \in \ft_0$ such that $(k t)_p = \ft_0$: then $t \in \Omega''$.

We will prove below that $\fn_x \subseteq \fg_{(0)}$ for all $x \in \Omega$. First, let us show how this property yields the result we want. Set $\mu = n\!-\!1$, and let $\beta$ be a $\mu$-linear alternating form on $\fg$, whose kernel contains $\fg_{(0)}$. Let $x \in \Omega$. For all $a_1,\ldots,a_{\mu} \in (kx)_p$ with Jordan-Chevalley decompositions $a_i = s_i + n_i$, we have $n_i \in \fn_x \subseteq \fg_{(0)}$. Hence, we obtain $\beta(a_1,\ldots,a_\mu) = \beta(s_1,\ldots,s_\mu)$. Now consider $\Omega_\beta = \{x�\in \Omega\ ;\ \beta(x,x^{[p]},\ldots,x^{[p]^{\mu-1}}) \neq 0\}$. This is again an open subset of $\fg$, and we can choose $\beta$ so that it is not empty. Using Equation (\ref{expressionwithbeta}), we obtain 
\begin{eqnarray*}
\varphi_i(x) & = &  \varphi_i(x_s) =  - \frac{\beta(x_s,\ldots,x_s^{[p]^\mu},\ldots,x_s^{[p]^{\mu-1}})^{p}}{\beta(x_s,\ldots,x_s^{[p]^i},\ldots,x_s^{[p]^{\mu-1}})^{p}} \\
& = &  - \frac{\beta(x,\ldots,x^{[p]^\mu},\ldots,x^{[p]^{\mu-1}})^{p}}{\beta(x,\ldots,x^{[p]^i},\ldots,x^{[p]^{\mu-1}})^{p}}.
\end{eqnarray*}
for all $x \in \Omega_\beta$. This shows that $\varphi_i^{1/p}$ is a rational function on $\Omega_\beta$, whence $\varphi_i^{1/p} \in k(\fg)$. Since its $p$-th power is polynomial and $k[\fg]$ is integrally closed, we obtain $\varphi_i^{1/p} \in k[\fg]$.

Let $x \in \Omega$, so that $\ft_x \cap \fg_{(0)}=(0)$. Consequently, $\ft_x$ is conjugate to the generic torus $\ft_0$ by Theorem \ref{AGT2}, and we may assume that $\ft_x = \ft_0$. 
For all $i \in \{1,\ldots,n\}$, let $\xi_i := x_i + 1 \in \fA_n$. Thus $\xi_i^p = 1$ and $\xi_1,\ldots,\xi_n$ are invertible generators of $\fA_n$. We use the multi-index notation $\xi^{a} = 
\xi_1^{a_1} \cdots \xi_n^{a_n}$ when $a = (a_1,\ldots,a_n) \in \ZZ^n$. The ``partial derivative'' $\partial_i$ satisfies $\partial_i(\xi_j) = \delta_{i,j}$ for all $i,j$. Now for $i \in 
\{�1,\ldots,n\!-\!1\}$, let $\theta_i := \xi_i \partial_i - \xi_n \partial_n$, so that $\ft_0 = \sum_{i = 1}^{n-1} k \theta_i$. We also set $\theta_n := \xi_n \partial_n \in W(n)$. It is straightforward to check that
\begin{equation*}
\theta_i(\xi^{a}) = (a_i - a_n) \xi^{a}, \quad \forall\, i \in \{1,\ldots,n\!-\!1\}.
\end{equation*}
Hence, the subalgebra $\fA_n^{\ft_0} := \{x \in \fA_n \ ; \ t(x)= 0 \ \ \forall \ t \in \ft_0\}$ of $\ft_0$-constants in $\fA_n$ is given by 
\begin{equation*}
\fA_n^{\ft_0} = k[\zeta], \ \ \text{with} \ \  \zeta = \xi_1 \cdots \xi_{n-1} \xi_n^{-1}.
\end{equation*}

Let us compute the centralizer of $\ft_0$ in $W(n)$. It is easy to check that $\{\theta_1,\ldots,\theta_n\} $ is a basis of the free $\fA_n$-module $W(n)$. Write a derivation $D \in W(n)$ as $D = \sum_{i = 1}^n f_n \theta_i$, where $f_i \in \fA_n$ for all $i$. Note that $[\ft_0,\theta_n] = (0)$, so that $0 = [t,D] = \sum_{i = 1}^{n} t(f_i) \theta_i$ for all $t \in \ft_0$. It 
readily follows that
\[C_{W(n)}(\ft_0) = \bigoplus_{i=1}^{n} k[\zeta] \theta_i.\]
If, in addition, $D \in (kx)_p \subseteq C_{W(n)}(\ft_0) \cap S(n)$, then
\begin{equation*}
0 = \Div(D) =  \sum_{i=1}^n \Div(f_i \theta_i) =  \sum_{i=1}^n f_i \Div(\theta_i) + \theta_i(f_i)  = f_n + \theta_n(f_n).
\end{equation*}
Now $f_n = f_n(\zeta)$ is a polynomial expression in $\zeta$. Denote formally $f_n'(\zeta) = d f_n / d \zeta$. Then $\theta_n(f_n) = f_n'(\zeta) \theta_n(\zeta) = -\zeta f_n'(\zeta)$. Consequently, we must have $f_n(\zeta) = \zeta f_n'(\zeta)$, whence $f_n \in k \zeta$. Thus, there exists $\lambda \in k$ such that
\begin{equation*}
D = \sum_{i = 1}^{n-1} f_i \theta_i + \lambda \zeta \theta_n.
\end{equation*}
Thanks to the proof of \cite[(2.2.3)]{JMB}, we get that the ring of constants $\fA_n^x \neq k$. Since $\fA_n^x =  \fA_n^{(kx)_p} \subseteq \fA_n^{\ft_0} = k[\zeta]$, there exists a non-constant polynomial $g \in k[\zeta]$ such that $D(g) = 0$. Hence,
\begin{equation*}
0 = \sum_{i = 0}^{n-1} f_i \theta_i(g) + \lambda \zeta \theta_n(g) = \lambda \zeta \theta_n(g) = -\lambda \zeta^2 g'(\zeta).
\end{equation*}
Since $\zeta$ is invertible and $g \not \in k$, this forces $\lambda = 0$. We have thus shown that
\begin{equation} \label{inclusionkxp}
(kx)_p \subseteq \sum_{i = 1}^{n-1} k[\zeta] \theta_i.
\end{equation}

The projection $\fg \twoheadrightarrow \fg / \fg_{(0)}$ induces a linear map $(kx)_p \to \fg / \fg_{(0)}$. Using the inclusion (\ref{inclusionkxp}) we see that its image is spanned by the images of $\theta_1,\ldots,\theta_{n-1}$, and hence has dimension $n\!-\!1$. Therefore the kernel $(kx)_p \cap \fg_{(0)}$ has dimension 1. On the other hand, $(kx)_p \cap \fg_{(0)}$ is a $p$-subalgebra of $(kx)_p = \ft_0 \oplus \fn_x$ which doesn't intersect $\ft_0$: so it must be $\fn_x$, and $\fn_x \subseteq \fg_{(0)}$, as we wanted to show. \end{proof}

\end{document}